\documentclass[11pt]{article}
\usepackage[hmargin=1.0in]{geometry}
\usepackage{amsmath,amsfonts,graphicx,subfigure,color,multicol}
\usepackage{empheq}
\usepackage{amssymb}
\usepackage{epstopdf}
\usepackage{float,topcapt}
\usepackage{hyperref}
\usepackage{multirow}
\usepackage{rotating}
\usepackage{longtable}
\usepackage{algorithm}
\usepackage{algpseudocode}
\usepackage{relsize}
\usepackage{pdflscape}
\usepackage{placeins}
\usepackage[font=small,labelfont=bf]{caption}
\usepackage{url}

\DeclareMathOperator*{\maximize}{maximize}
\DeclareMathOperator*{\minimize}{minimize}

\begin{document}

\newtheorem{problem}{Problem}

\newtheorem{remark}{Remark}

\newtheorem{dfn}{Definition}

\floatstyle{boxed}
\newfloat{algorithm}{h}{alg}
\floatname{algorithm}{Algorithm}

\newcommand{\nablav}{\vec{\nabla}}
\newcommand{\xv}{\vec{x}}
\newcommand{\qb}{\mathbf{q}}
\newcommand{\sbb}{\mathbf{s}}
\newcommand{\fbv}{\vec{\textbf{f}}}
\newcommand{\ximap}{\vec{\xi}}
\newcommand{\qbmap}{\qb^{\xiexp}}
\newcommand{\xiexp}{\hspace{0.08cm}\ximap}
\newcommand{\fbvmap}{\fbv^{\xiexp}}
\newcommand{\Jmap}{\vec{\vec{\left.\mathrm{J}\right.}}}
\newcommand{\nablavmap}{\vec{\nabla}^{\xiexp}}

\newcommand{\Qb}{\mathbf{Q}}
\newcommand{\Fbv}{\vec{\textbf{F}}}
\newcommand{\Fbvmap}{\vec{\textbf{F}}^{\xiexp}}
\newcommand{\Fbvnum}{\vec{\textbf{F}}_\mathrm{num}^{\xiexp}}
\newcommand{\nmap}{\vec{n}^{\xiexp}}

\newcommand{\rhob}{{\rho}_{0}}
\newcommand{\pb}{{p}_{0}}
\newcommand{\cb}{{c}_{0}}

\newcommand{\Real}{\mathbb{R}}
\newcommand{\Complex}{\mathbb{C}}

\newcommand{\yb}{\mathbf{y}}
\newcommand{\eb}{\mathbf{e}}
\newcommand{\Ab}{\mathbf{A}}
\newcommand{\bb}{\mathbf{b}}

\newcommand{\ub}{\mathbf{u}}
\newcommand{\Fb}{\mathbf{F}}
\newcommand{\Lb}{\mathbf{L}}

\newcommand{\dt}{\Delta t}
\newcommand{\dtm}{\nu_{\mathrm{stab}}}
\newcommand{\dta}{\nu_{\mathrm{acc}}}
\newcommand{\dr}{\Delta r}

\newcommand{\abs}[1]{\vert #1 \vert}
\newcommand{\maxnorm}[1]{\vert \vert #1 \vert \vert_{\infty}}

\newcommand{\ndof}{N^\text{DOF}}
\newcommand{\Oop}{\mathcal O}

\newcommand{\chistab}{\chi_{\mathrm{stab}}}
\newcommand{\chiacc}{\chi_{\mathrm{acc}}}

\newcommand{\e}[1]{\ensuremath{\times 10^{#1}}}

\title{Optimized explicit Runge-Kutta schemes for the spectral difference method
applied to wave propagation problems}
\author{M. Parsani\thanks{King Abdullah University of Science and Technology (KAUST),
    Division of Mathematical and Computer Sciences and Engineering,
    Thuwal, 23955-6900. Saudi Arabia  (\mbox{matteo.parsani@kaust.edu.sa}).} \and
    David I. Ketcheson\thanks{King Abdullah University of Science and Technology (KAUST),
    Division of Mathematical and Computer Sciences and Engineering,
    Thuwal, 23955-6900. Saudi Arabia.} \and
    W. Deconinck\thanks{Vrije Universiteit Brussel, Department of Mechanical Engineering,
    Pleinlaan 2, 1050 Brussels. Belgium.}}

\maketitle

\begin{abstract}Explicit Runge--Kutta schemes with large stable step sizes are
developed for integration of high order spectral difference
 spatial discretizations on quadrilateral grids.
The new schemes permit an effective 
time step that is substantially larger than the maximum admissible 
time step of standard explicit Runge--Kutta schemes available in literature.
Furthermore, they have a small principal error norm and
admit a low-storage implementation. 
The advantages of the new schemes are demonstrated through application to 
the Euler equations and the linearized Euler equations.
\end{abstract}

\section{Introduction} \label{sec:intro}
Throughout the past two decades, the development of high-order accurate spatial 
discretization has been one of the major fields of research in 
numerical analysis, computational
fluid dynamics (CFD), computational aeroacoustics (CAA), computational 
electromagnetism (CEM) and in general computational physics
characterized by linear and nonlinear wave propagation phenomena. 
High-order discretizations have the potential to improve the computational 
efficiency required to achieve a desired error level by allowing the use of coarser 
grids. Indeed, in modern wave propagation problems characterized by complicated 
geometries, complex physics and a wide disparity of length scales (e.g., large eddy 
simulation, turbulent combustion, flow around flapping wings, rotor-blade 
interaction), the need for high-accuracy solutions 
leads to a prohibitive computational cost when low-order (i.e., first- and 
second-order accurate) discretizations are used. High-order 
schemes have much better wave propagation properties and a truncation error that
decreases more rapidly than that of low-order schemes if the solution is 
sufficiently smooth. Therefore, for problems that require very low numerical 
dissipation and small error levels, it may be advantageous to use high-order 
spatial discretization schemes; see for instance 
\cite{Chen2005b,Colonius2004,Wang2007e,Vincent2011}.

Among high-order methods, the spectral difference (SD) scheme 
is receiving increasing attention \cite{Cox2012,Gao2012,Liang2102,Ou2011,Parsania,Parsani2010a,
Balan20122359,Lodato2012}. 
The SD scheme offers several interesting properties. It is able to obtain solutions
with arbitrarily high order of accuracy. It can be applied to unstructured quadrilateral
and hexahedral meshes (tensorial cells). The conservation laws to be solved are in differential form,
avoiding the use of costly high-order accurate quadrature formulas. 

Although the formulation of high-order spatial discretization is now fairly 
mature, the development of techniques for efficiently solving systems of 
ordinary differential equations (ODEs) arising from high-order accurate spatial
discretizations has received less attention. The cost of solving an initial 
value problem up to a fixed time is inversely
proportional to the time step used, so it is desirable to use the largest
step size possible if the temporal discretization error is
acceptable. For higher order schemes, the spectrum of the Jacobian of
the semi-discretization often has increasingly large eigenvalues.  As a result,
the step size is often limited by stability requirements, which become stricter 
with higher order methods.
Implicit methods allow the use of much larger step sizes, but lead to very large memory 
requirements that may not be feasible. 
Moreover, the development of efficient algebraic
solvers for high-order implicit discretizations remains challenging. 
Thus, explicit time integration methods that allow large step sizes and require
less memory seem to be an appealing alternative.

This work focuses on the development of new optimized explicit Runge-Kutta (ERK) 
schemes to compute wave propagation efficiently and accurately with 
high-order SD methods on unstructured uniform or quasi-uniform quadrilateral cell grids. 
The schemes are optimized with respect to the spectrum of the SD discretization,
using the two-dimensional (2D) advection equation as a model problem.
Linear stability optimization determines the coefficients of the stability
polynomial but does not uniquely determine the full RK method.
A second optimization step is used to determine the 
Butcher coefficients of the scheme, optimized for a small leading truncation
error constant and low-storage form. The low-storage form is crucial for memory
reasons, since many stages are used.

Many authors have studied the design of optimal ERK schemes with many stages for integration
of high order discretizations of partial differential equations (PDEs).  
Past efforts focused on schemes with a relatively smaller number of stages
\cite{Allampalli2009,toulorge2011b,diehl2010comparison,Niegemann2011,Bernardini2009b}.
By using the algorithm developed in \cite{Ketcheson2012a}, we are able to develop
schemes with much larger number of stages and with higher order of accuracy.
Our work is also the first to develop schemes specifically for the spectral difference
semi-discretization.
Whereas past studies have focused on application to linear problems only, and typically
employed structured grids, we validate the effectiveness of our methods also
on a nonlinear, fully unstructured example.
Our new optimal ERK methods increase the computational efficiency of the
SD method for wave propagation problems up to $65\%$ and $57\%$, respectively
for 4$^\mathrm{th}$-order and 5$^\mathrm{th}$-order accurate simulations.

The remainder of the paper is organized as follows. 
In Section \ref{sec:sd-discretization}, we review the SD method for tensor product
cells (quadrilateral and hexahedral cells).  Section \ref{sec:rk-opt} is 
devoted to the description of the two-step optimization algorithm used to 
design new ERK schemes, in which we first select an optimal stability polynomial
and then design a corresponding ERK method. 
We also present the main features of the optimized methods
and discuss their theoretical efficiency. Section \ref{sec:applications}
presents numerical results for three benchmark test problems, which
demonstrate that the new schemes lead to large performance gains
over standard ERK schemes available in the literature,
for both linear and nonlinear problems.
Conclusions and future directions are given in Section \ref{sec:conclusions}.

\section{Spectral difference discretization}\label{sec:sd-discretization}
In this section, we review the spectral difference approach to semi-discretization
of hyperbolic conservation laws.

Consider the general hyperbolic system of conservation laws over a 
three-dimensional domain $\Omega \subset\mathbb{R}^{3}$ with boundary 
$\partial\Omega$ and completed with consistent initial and boundary conditions:
\begin{equation} 
    \begin{cases}
\frac{\partial \qb}{\partial t} + \nablav\cdot \fbv\left(\qb\right) = \sbb(\qb) & \textrm{ in } \Omega \times [t^0,t^e] \\ 
\qb\left(\xv,0\right) = \qb^0\left(\xv\right) & \textrm{ on } \Omega \\
\qb\arrowvert_{\partial \Omega}\left(t\right) = \qb^b\left(t\right) & \textrm{ on } \partial\Omega, 
    \end{cases}
\label{eq:genHypSys}
\end{equation}
Here, $\xv$, $\qb$, $\fbv$, $\sbb$, $t^0$ and $t^e$ are respectively the position vector, 
the vector of the conserved variables, the flux vector, source terms, and the lower and upper 
bound of the time interval. The spatial domain $\Omega$ is discretized into
tensor product cells with domain and 
boundary $\Omega_i$ and $\partial \Omega_i$. 

For each cell $i$, take a mapped coordinate system 
$\vec{\xi} = \left[\xi,\eta,\zeta\right]^T$. The transformation from the 
standard to the physical element in the global Cartesian coordinates for the 
cell $i$ is given by 
\begin{equation}
    \xv_i= \left[\begin{array}{c} x_i\left(\xi,\eta,\zeta\right)\\
    y_i\left(\xi,\eta,\zeta\right) \\ z_i\left(\xi,\eta,\zeta\right)\\
    \end{array}\right]
    =\xv_i\left(\ximap\right),
\label{eq:mapcoordtrans}
\end{equation}
with Jacobian matrix $\Jmap_i$ and Jacobian determinant 
$J_i$. The fluxes projected in the mapped coordinate system 
($\fbvmap_i$) are then related to the flux components in the global coordinate 
system by
\begin{equation}
\fbvmap_i = 
\left[\begin{array}{l} 
\mathbf{f}_i^{\xiexp} \\
\mathbf{g}_i^{\xiexp} \\ 
\mathbf{h}_i^{\xiexp}
\end{array}\right]
= J_i\Jmap_i^{-1}
\left[\begin{array}{l}
\mathbf{f}_i\\
\mathbf{g}_i\\
\mathbf{h}_i
\end{array}\right]
= J_i\Jmap_i^{-1} \fbv_i.
\label{eq: fluxMap}
\end{equation}
Therefore, the hyperbolic system \eqref{eq:genHypSys} can be written in 
the mapped coordinate system as
\begin{equation}
\frac{\partial \left(J_i\qb\right)}{\partial t} \equiv \frac{\partial\qb_i^{\xiexp}}{\partial t} 
= -\frac{\partial\mathbf{f}_i^{\xiexp}}{\partial\xi}
-\frac{\partial\mathbf{g}_i^{\xiexp}}{\partial\eta}
-\frac{\partial\mathbf{h}_i^{\xiexp}}{\partial\zeta} 
= -\nablavmap\cdot\fbvmap_i,
\label{eq:genHypSysMap}
\end{equation}
where $\qbmap_i \equiv J_i \qb$ and $\nablavmap$ are the conserved variables and
the differential operator in the mapped coordinate system, respectively.

For a $\left(p+1\right)$-th-order accurate $d$-dimensional scheme, 
$N^s$ \emph{solution collocation points} with index $j$ are introduced at positions 
$\ximap_j^s$ in each cell $i$, with $N^s$ given by $N^s=\left(p+1\right)^{3}$.
Given the values at these points, a polynomial approximation of degree $p$ of
the solution in cell $i$ can be constructed. This polynomial is called the
\emph{solution polynomial} and is usually composed of a set of Lagrangian 
basis polynomial $L_j^s\left(\ximap\right)$ of degree $p$:
\begin{equation}
     \begin{array}{ll}
	    \Qb_i\left(\ximap\right) &= \sum_{j=1}^{N^s} \Qb_{i,j} \, L_j^s\left(\ximap\right) \\
        L_j^s\left(\ximap^s\right) &= \delta_{jm}, \quad j,m=1,...,N^s.
     \end{array}
\label{eq:sdsolpoly}
\end{equation}
Therefore, the interpolation coefficients are given as 
$\Qb_{i,j} = \Qb_i\left(\ximap_j^s\right)$ where $\Qb_{i,j}$ are the conserved 
variables at the solution points, i.e. the unknowns of the SD method. 

The divergence of the mapped fluxes 
$\nablavmap\cdot\fbvmap$ at the solution points is computed by 
introducing a set of $N^f$ flux collocation points with index $l$ and at positions 
$\ximap_l^f$, supporting a polynomial of degree $p+1$. The evolution of the 
mapped flux vector $\fbvmap$ in cell $i$ 
is then approximated by a flux polynomial $\Fbvmap_i$, which is obtained
by reconstructing the solution variables at the flux points and evaluating the 
fluxes $\Fbvmap_{i,l}$ at these points. The flux is represented by a separate 
Lagrange polynomial:
\begin{equation}
     \begin{array}{ll}
        \Fbvmap_i\left(\ximap\right) = \sum_{l=1}^{N^f} \Fbvmap_{i,l} \, L_l^f\left(\ximap\right) \\
        L_l^f\left(\ximap_m^f\right) = \delta_{lm}, \quad l,m=1,...,N^f.
    \end{array}
\label{eq:sdfluxpoly}
\end{equation}
Hence, the coefficients of the flux interpolation are defined as
\begin{equation}
    \Fbvmap_{i,l} =
        \begin{cases}
            \Fbvmap_i\left(\ximap_l^f\right), & \quad \ximap_l^f \in \Omega_i \\
            \Fbvnum\left(\ximap_l^f\right), & \quad \ximap_l^f \in \partial\Omega_i
        \end{cases}
\end{equation}
where $\Fbvnum$ is the flux vector at the cell interface. In order to maintain 
conservation at a cell level, the flux component normal to 
a face (i.e. $\Fbvnum \cdot \nmap$) must be continuous between two neighboring 
cells. However, the solution
at a face is in general not continuous and requires the solution of a Riemann 
problem. For many nonlinear hyperbolic systems, such as the 
compressible Euler equations, the exact Riemann solution cannot be
written in closed form and is prohibitively expensive to compute. Therefore,
cheaper approximate Riemann solvers are typically used. The tangential component
of $\Fbvnum$ is usually taken from the interior cell (see for instance 
\cite{Wang2006a}).

Taking the divergence of the flux polynomial $\nablavmap\cdot\Fbvmap_i$ in the 
solution points results in the following modified form of 
\eqref{eq:genHypSysMap}, describing the evolution of the conservative
variables in the solution points:
\begin{equation}
\frac{d\Qb_{i,j}}{dt} = 
-\left.\nablav\cdot\Fbv_i\right|_j = -\frac{1}{J_{i,j}}\left.\nablavmap\cdot
\Fbvmap_i\right|_j = \mathbf{R}_{i,j},
\label{eq:residualsSD}
\end{equation}
where $\Fbv_i$ is the flux polynomial vector in the physical space whereas 
$\mathbf{R}_{i,j}$ is the SD residual associated with $\Qb_{i,j}$. This is
a system of ODEs, in time, for the unknowns 
$\Qb_{i,j}$.

\subsection{Solution and flux points distributions} \label{subsec:solFluxPoints}
Huynh \cite{THuynh2007} showed
that for quadrilateral and hexahedral cells, tensor product flux point 
distributions based on a one-dimensional (1D) flux point distribution consisting of the end points
and the Legendre-Gauss quadrature points lead to stable schemes for arbitrary
order of accuracy.

In 2008, Van den Abeele et al.~\cite{Abeele2008a} showed an interesting
property of the SD method, namely that it is independent of the positions 
of its solution points in most general circumstances, for both simplex and 
tensor-product cells. In the above work it has been shown that the distribution 
of the solution points has very little influence on the properties of the SD 
schemes, and in fact, for linear problems, different distributions lead to 
identical results. This property greatly simplifies the design of SD schemes, 
since only the flux point distributions has to be taken care with. It also implies 
an important improvement in efficiency, since the solution points can be placed 
at flux point positions and thus a significant number of solution reconstructions can be 
avoided. Recently, this property has been proved by Jameson 
\cite{Jameson2010}.

\begin{figure}[ht]
\centering
\includegraphics[scale=0.3]{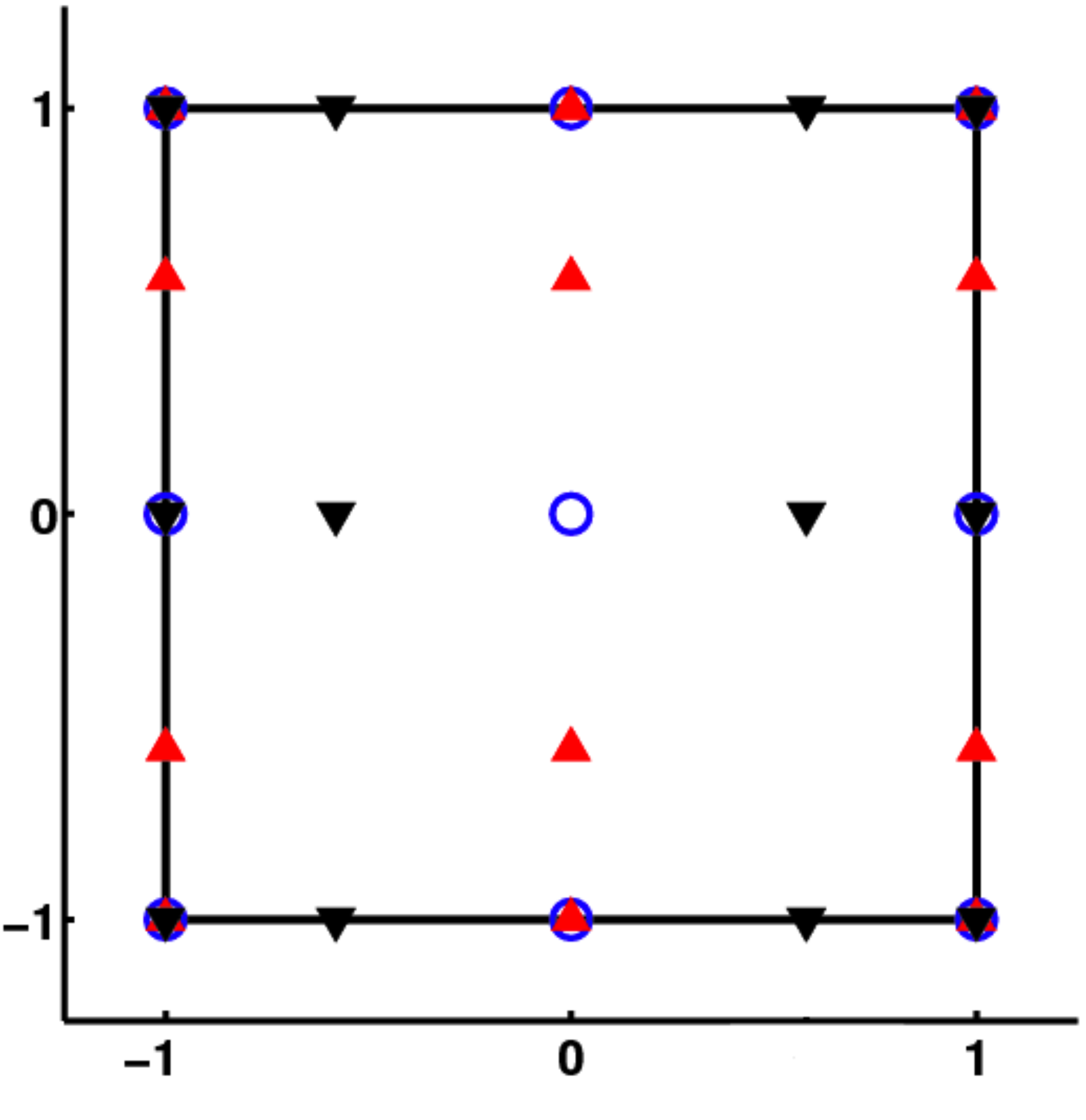}
\caption{Typical 3$^\mathrm{rd}$-order ($p = 2$) quadrilateral SD cells, 
with component-wise flux point distributions. Solution 
(\textcolor{blue}{$\circ$}), $\xi$- ($\blacktriangledown$) and 
$\eta$-flux points ({\color{red}{$\blacktriangle$}}).}
\label{fig: thirdordersd2pract}
\end{figure}

Figure \ref{fig: thirdordersd2pract} shows a typical distribution of flux and solution points for a 
third-order SD scheme in 2D.

\subsection{Advection equation} \label{sec:advection}
The advection equation represents the simplest hyperbolic conservation law.
It models the advection of a scalar conserved variable $q$ with constant advection speed 
$\vec{a}$. The conserved variables and the convective flux are then
\begin{equation}
\begin{aligned}
  \mathbf{q} &= q, \\
  \vec{\mathbf{f}} & = \vec{a}\,q.
\label{eq: linadv}
\end{aligned}
\end{equation}
Therefore, the conservation law reads
\begin{equation}
	\frac{\partial q}{\partial t} + \vec{\nabla} \cdot(\vec{a}\,q) = 0.
\label{eq:lin-adv}
\end{equation}

\subsubsection{Spectrum of the two-dimensional advection equation}\label{subsub:spectrum}
Despite its simplicity, the one-dimensional advection equation is often
used as a model problem to design and analyze new
spatial discretization schemes for convection-dominated problems. It is also used to optimize the 
coefficients of time integration algorithms for CFD, and for wave propagation 
problems in general.  

In this work, however, we use the 2D advection equation as a model.
In the 2D case, the discrete operator
arising from the spatial discretization is also a function of the convective 
velocity direction. This approach allows us to consider different wave propagation 
trajectories and optimize the RK coefficients by using a richer 
spectrum (or Fourier footprint) than that of the 1D advection 
equation. The richer spectrum leads to a design of more robust schemes.
In the remaining part of this section, the procedure used to compute the
Fourier footprint is described.

\begin{figure}[ht!]
\centering
\includegraphics[scale=0.5]{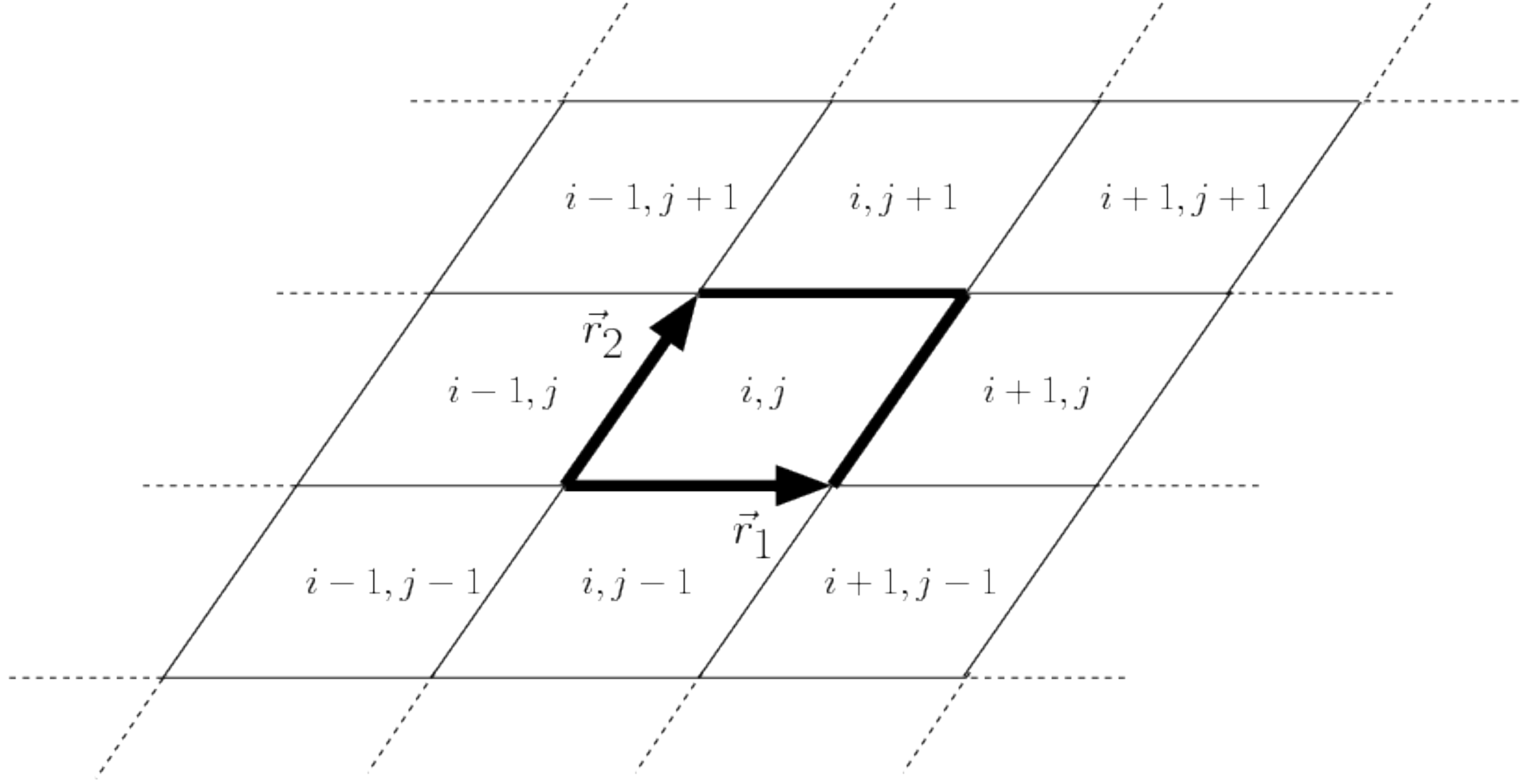}
\caption{Generating pattern.}
\label{fig:gen-patt-vn}
\end{figure}

Equation (\ref{eq:lin-adv}) is discretized in space by the SD scheme.
A uniform grid with periodic boundary conditions is considered. The
grid is defined by a \emph{generating pattern}, which is the smallest part from
which the full grid can be reconstructed by periodically repeating the pattern 
in all directions. For the 2D case and uniform 
quadrilateral meshes, the generating pattern is completely defined by the vectors $\vec{r}_1$ and 
$\vec{r}_2$ (see Figure~\ref{fig:gen-patt-vn}) whose non-dimensional form 
is obtained by scaling them with 
the length of $\vec{r}_1$, denoted by 
$\Delta r$: $\vec{r}_1 \equiv \Delta r \, \vec{r}_1^{\,\prime}$ and $\vec{r}_2 \equiv \Delta r \, \vec{r}_2^{\, \prime}$. 
If the dimensionless vector $\vec{r}_1^{\,\prime}$ is chosen as $\left[1 \,\, 0\right]^T$, then
the dimensionless mesh is completely defined by the two components of $\vec{r}_2^{\,\prime }$.

The advection speed $\vec{a}$ in Equation (\ref{eq:lin-adv}) is 
defined by its amplitude $|\vec{a}|$ and orientation angle $\psi$:
\begin{equation}\label{eq:advection-velocity-def}
  \vec{a} = |\vec{a}|\left[\begin{array}{c}\cos\psi \\ \sin\psi\end{array}\right].
\end{equation}

At cell faces the solution is discontinuous, so two values for the convected
variables are available.
The normal flux component is calculated using the following approximate Riemann solver:
\begin{eqnarray}
\vec{F}\left(Q_L,Q_R\right) \cdot \vec{1}_n = \vec{a} \cdot \vec{1}_n\frac{Q_L + Q_R}{2} - \left|\vec{a} \cdot \vec{1}_n\right| \frac{Q_R - Q_L}{2},
\label{eq:num-riemann-flux}
\end{eqnarray}
where $\vec{1}_n$ is the unit normal oriented from the left to the right side 
and indices $L$ and $R$ indicate respectively the left and right neighboring 
cell to a face. In the present 
analysis, the internal component in each cell is used for the tangential flux 
component.

After the SD semi-discretization of (\ref{eq:lin-adv}) on a 
uniform quadrilateral mesh, the following system of ODEs is obtained:
\begin{align}
\frac{d \mathbf{Q}_{i,j}}{d t} + \frac{|\vec{a}|}{\Delta r} \left( \mathbf{T}^{0,0} \, \mathbf{Q}_{i,j} \right. &+ \mathbf{T}^{-1,0} \, \mathbf{Q}_{i-1,j} + \mathbf{T}^{0,-1} \,  \mathbf{Q}_{i,j-1} \nonumber \\
& \left.+ \mathbf{T}^{+1,0} \, \mathbf{Q}_{i+1,j} + \mathbf{T}^{0,+1} \, \mathbf{Q}_{i,j+1} \right)= 0,
\label{eq:lin-adv-vn}
\end{align}
where the five matrices $\mathbf{T}$ are determined by the coefficients of the spatial 
discretization. They depend on the order of accuracy $p$ of the SD scheme, the generating pattern, and the advection 
velocity orientation angle $\psi$. The column vector 
$\mathbf{Q}_{i,j}$ contains all solution point variables of the cell with indices $i$ and $j$ (Figure~\ref{fig:gen-patt-vn}).

Inserting the following plane Fourier wave
\begin{eqnarray}
	\mathbf{Q}_{i,j}\left(t\right) &=& \tilde{\mathbf{Q}}\left(t\right) e^{I\,\vec{k}\cdot\left( i\,\vec{r}_1^{\,\prime} + j\,\vec{r}_2^{\,\prime} \right)\Delta r} \nonumber\\
&=& \tilde{\mathbf{Q}}\left(t\right) e^{I\,\vec{K}\cdot\left( i\,\vec{r}_1^{\,\prime} + j\,\vec{r}_2^{\,\prime} \right)}
\label{eq:fourier-wave-vn}
\end{eqnarray}
into Equation (\ref{eq:lin-adv-vn}) results in
\begin{align}
\frac{d\tilde{\mathbf{Q}}}{dt} + \frac{|\vec{a}|}{\Delta r} \left(\mathbf{T}^{0,0} \right. &+ \mathbf{T}^{-1,0} \, e^{-I\vec{K}\cdot \vec{r}_1^{\,\prime}} + \mathbf{T}^{0,-1} \, e^{-I\vec{K}\cdot \vec{r}_2^{\,\prime}} \nonumber \\
& \left. + \mathbf{T}^{+1,0} \, e^{I\vec{K}\cdot \vec{r}_1^{\,\prime}} + \mathbf{T}^{0,+1} \, e^{I\vec{K}\cdot \vec{r}_2^{\,\prime}} \right) \tilde{\mathbf{Q}} = 0,
\label{eq:odes-vn}
\end{align}
where $I\equiv\sqrt{-1}$ is the imaginary unit number. Here, $\vec{k}$ and $\vec{K}$ 
are the wave vector and the dimensionless wave vector given by
\begin{equation}
	\vec{k}=|\vec{k}|\,\left[\begin{array}{c}\cos\theta \\ \sin\theta\end{array}\right]
\end{equation}
and 
\begin{equation}\label{eq:dimensionless-K}
	\vec{K}=\vec{k}\,\Delta r = |\vec{k}|\,\Delta r \left[\begin{array}{c}\cos\theta \\ \sin\theta\end{array}\right], 
\end{equation}
where $\theta$ is the angle between the wave vector $\vec{k}$ and the horizontal 
axis.

Equation (\ref{eq:odes-vn}) can be written as 
\begin{equation}
\frac{d\tilde{\mathbf{Q}}}{dt} = \frac{|\vec{a}|}{\dr}\ \Lb\ \tilde{\mathbf{Q}},
\label{eq:odes-compact-vn}
\end{equation}
where the matrix $\Lb$ is defined by the SD spatial operator. The set of eigenvalues of 
$\Lb$ is the spectrum or Fourier footprint of the spatial
discretization. The spectrum depends 
on the order of accuracy $p$ of the SD scheme, the 
generating pattern, the direction $\psi$ of the convective velocity, and the 
dimensionless wave number vector $\vec{K}$.

Here we take a uniform Cartesian grid defined by
\begin{equation}
\vec{r}_{1}^{\, \prime }=\left( 
\begin{array}{c}
\dr \\ 
0
\end{array}
\right), \hspace{1cm}
\vec{r}_{2}^{\, \prime }=\left( 
\begin{array}{c}
0 \\ 
\dr
\end{array}
\right)
\end{equation}
and, for a given order of accuracy, we compute the spectrum of the operator $\Lb$ 
by varying $\psi$, $K$ and $\theta$. 

\section{Optimized Runge--Kutta schemes}\label{sec:rk-opt}
The spectral difference semi-discretization of a PDE described in the previous 
section leads to an initial value problem
\begin{equation} \label{eq:ivp}
  \left\{
  \begin{split}
    \Qb'(t) & = \Fb(\Qb)  \\
    \Qb(0) & = \Qb_0
  \end{split}
  \right.,
\end{equation}
where $\Qb(t) : \Real\to\Real^{\ndof}$ and $\Fb: \Real^{\ndof} \to 
\Real^{\ndof}$ are the vector of the unknowns and the vector of the residuals,
respectively. The number of degrees of freedom
is denoted by $\ndof = N\times N^{s}$, where $N$ is the number of cells 
used to discretize the domain $\Omega$.
System \eqref{eq:ivp} is typically integrated by using a high-order accurate ERK
time discretization, which takes the form
\begin{equation}
\begin{aligned} \label{eq:rk}
    Y_i & = Q^n + \dt \sum_{j=1}^{i-1} a_{ij} F(Y_j) \\
    Q^{n+1} & = Q^n + \dt \sum_{j=1}^{i-1} b_{j} F(Y_j),
\end{aligned}
\end{equation}
for a scalar ODE.
The properties of the RK method are determined by its coefficient
matrix $\Ab = [a_{ij}]$ and column vector $\bb = [b_j]$
which are referred to as 
the Butcher coefficients \cite{Butcher2008}.
In this section, we describe our approach to designing RK schemes that
maximize the absolutely stable time step size, have reasonably small error
constants, and can be implemented with low storage requirements.

\subsection{Optimization of the stability polynomial}\label{subsec:optim-lin-stab}
Stability of RK integration is studied by applying the method \eqref{eq:rk}
to the linear scalar test problem $Q'(t)=\lambda\ Q$.
Any RK method applied to this problem yields an iteration of the form
\begin{align} \label{eq:RKiter}
Q^{n+1} = \psi(\dt\ \lambda)\ Q^{n}
\end{align}
where the {\em stability function} $\psi(z)$ depends only 
on the coefficients of the
RK method (\cite[Section 4.3]{Gottlieb2011}\cite{Butcher2008,Hairer1991}):
\begin{align} \label{eq:psi}
\psi(z) = 1 + \sum_{j=0}^s \bb^T \Ab^{j-1} \eb \, z^j.
\end{align}
Here $\eb$ is a column vector of size $s$ made by ones.
The stability function governs the local
propagation of errors, since any perturbation to the solution will be multiplied
by $\psi(\dt\ \lambda)$ at each subsequent step.

We say the iteration \eqref{eq:RKiter} is absolutely stable if 
\begin{align} \label{eq:spec-cond}
\dt\ \lambda & \in S  & \text{where } & & S = \{ z \in \Complex : |\psi(z)|\le 1 \}.
\end{align}
The set $S$ is referred to as the {\em absolute stability region}.

When applied to a linear system of PDEs (such as the advection or linearized Euler
equations discussed in the previous section), the SD semi-discretization
leads to a linear, constant-coefficient initial value problem; i.e. \eqref{eq:ivp} with
$\Fb(\Qb) = \frac{|\vec{a}|}{\dr}\ \Lb\ \Qb$ where $\Lb$ is a fixed square matrix that
depends on the order of accuracy $p$ of the spatial discretization, the 
generating pattern (see Figure \ref{fig:gen-patt-vn}), the direction $\psi$ of 
the convective velocity defined in Equation \eqref{eq:advection-velocity-def}, 
and the dimensionless wave number vector $\vec{K}$ given by Equation
\eqref{eq:dimensionless-K}.

Application of a RK method to \eqref{eq:ivp}
leads to the iteration 
\begin{align} \label{eq:RKiter2}
\Qb^{n+1} = \psi(\nu\Lb)\ \Qb^{n}
\end{align}
where $\nu=|\vec{a}|\,\frac{\dt}{\dr}$ is the CFL number.
Assume that $\Lb$ is diagonalizable and let $\lambda_i, i=1,\dots,\ndof$ denote its eigenvalues.
Then the solution is absolutely stable for CFL number $\nu$ if
\begin{align} \label{eq:spec-cond-sys}
\nu\lambda_i \in S & \ \ \mbox{ for } 1\le i \le \ndof.
\end{align}
Thus the maximum absolutely stable step size is
\begin{align}
\dtm = \max \{ \nu\ge0 : |\psi(\nu\lambda_i)|\le 1 \mbox{ for } i=1,\dots,\ndof\}.
\end{align}
Although this analysis is based on the linear problem, 
it is often used to obtain a practical step size restriction for
the nonlinear problem \eqref{eq:ivp} by considering the spectrum of the Jacobian 
of $\Fb(\Qb)$.

In general, $\psi(z)$ for an $s$-stage, order $p$, ERK method is
a polynomial of degree $s$ that differs from the exponential function by terms of
order $z^{p+1}$:
\begin{align} \label{eq:polyform}
\psi(z) & = \sum_{j=0}^s \beta_j z^j = \sum_{j=0}^p \frac{1}{j!} z^j 
    + \sum_{j=p+1}^s \beta_j z^j.
\end{align}
Comparing \eqref{eq:polyform} with \eqref{eq:psi}, we see that $\beta_j=\bb^T\, \Ab^{j-1}\,\eb$.

It is natural then to design optimal polynomials by choosing the coefficients $\beta_j$ in
\eqref{eq:polyform} so as to maximize $\dtm$.
The optimization problem may be stated formally as follows
\begin{problem}[Stability polynomial optimization] \label{prob:stabpoly}
    \begin{equation*}
      \begin{aligned}
        \textup{Choose } \{\beta_{p+1},\dots,\beta_s\} \textup{ to } & \maximize && \nu \\
        & \textnormal{subject to} && \\
        & && |\psi(\nu\lambda)|\le 1 \ \ \ \mbox{for all } \lambda \in\sigma(\Lb) \\
        & &&  \psi(z)-\exp(z) = \Oop(z^{p+1}).
      \end{aligned}
    \end{equation*}
\end{problem}
We solve Problem \ref{prob:stabpoly} using a convex optimization approach and 
bisection with respect
to the CFL number $\nu$, as described in \cite{Ketcheson2012a}. 
Specifically, we fix the step size $\nu$ and solve the resulting convex
feasibility problem, to determine whether there exists a choice of coefficients
$\beta$ that satisfy the constraints.  The upper bound for the initial bisection
interval is $10\,s$, and we use a monomial basis.  The bisection search is carried 
out to an accuracy of $10^{-7}$.

This approach allows us to optimize methods with large numbers of stages in
order to improve the maximum absolutely stable time step $\dtm$.
The optimization is carried out for 2$^\mathrm{nd}$- to 5$^\mathrm{th}$-order
accurate schemes;
the constraint points $\lambda$ are taken as the spectrum of the SD
semi-discretization of the same order of accuracy. The latter choice seems to be 
the most natural one. However, since in many 
situations the spatial error dominates and a temporal discretization should be
chosen to achieve a pre-defined error tolerance 
at the lowest possible cost, one might think that a very stable low order ERK 
schemes could be a valid alternative. As shown in 
\cite{Parsani-low-order-erk-sd}, by using a very low order of accuracy the linearly 
stable step size can be dramatically increased over the standard 
ERK schemes available in literature but such a gain reduces to about 5\% over the
optimal ERK methods presented in this work. Moreover, the speed-up is
obtained with some relatively large sacrifices in accuracy.

\begin{figure}[htbp!]
\centering
\subfigure[ERK(4,4).]{
\includegraphics[scale=0.3]{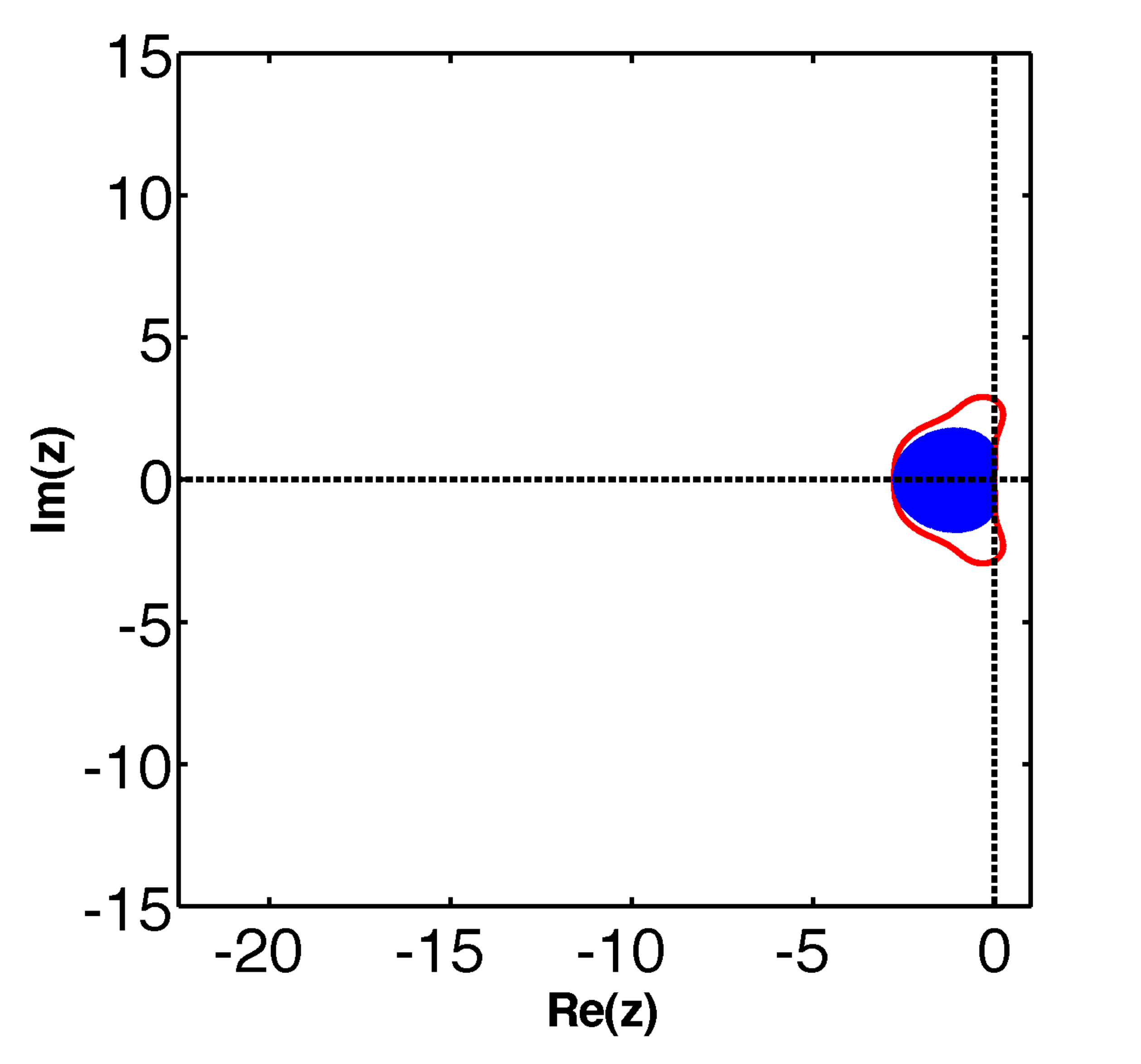}
\label{fig:scaled_spectrum_erk-4-4}}
\subfigure[ERK(18,4).]{
\includegraphics[scale=0.3]{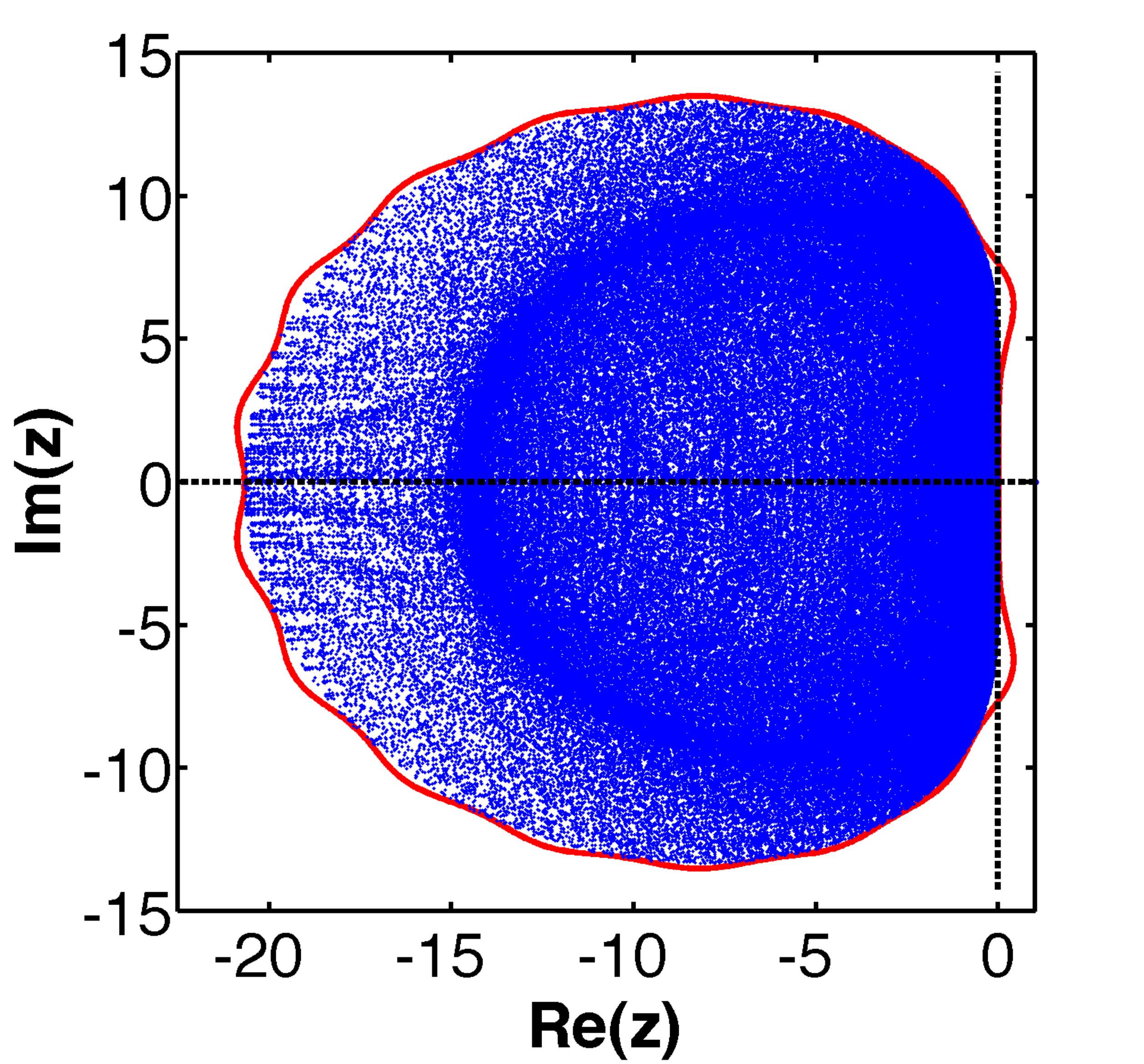}
\label{fig:scaled_spectrum_erk-18-4}}

\caption{Examples of stability region for 4$^\mathrm{th}$-order ERK methods (red lines)
and scaled Fourier footprint of the 4$^\mathrm{th}$-order SD scheme (blue dots);
$\psi\in[0,2\pi]$, $\theta \in[0,2\pi]$ and $|\vec{k}|\,\Delta r\in[0,2\pi]$.
\label{fig:example-stability-optimization}}
\end{figure}

Figure \ref{fig:example-stability-optimization} shows the stability regions 
of the classical 4-stage 4$^\mathrm{th}$-order  
(ERK(4,4)) \cite{Ku1901} and optimized 18-stage 4$^\mathrm{th}$-order 
(ERK(18,4)) methods superimposed on the Fourier footprint of the 
4$^\mathrm{th}$-order SD methods computed varying the direction $\psi$
of the wave propagation velocity vector $\vec{a}$,
the solution orientation defined by the angle $\theta$ and the dimensionless 
wave vector module $|\vec{k}|\,\Delta r$ (see Section \ref{subsub:spectrum}).
Clearly, the optimized method allows the use of a 
much larger step size.  Notice also that the stability region of the optimized
method has nearly the same shape as the convex hull of the SD spectrum.

\subsection{Determination of Runge--Kutta coefficients}\label{subsec:rk-coeffs}
The choice of stability polynomial does not fully determine the method; an ERK
method of $s$ stages has $s(s+1)/2$ coefficients and only $s$ of
them are constrained by the stability polynomial.
We now consider the problem of finding the RK coefficients $\Ab,\bb$
corresponding to a set of prescribed stability polynomial coefficients $\beta_j$.
We use the remaining degrees
of freedom to satisfy additional nonlinear order conditions, to obtain a
low-storage implementation, and ensure that the
truncation error coefficients are not too large.  

While the linear accuracy of the method is determined by the stability
polynomial, the nonlinear accuracy depends on larger set of order conditions
$\tau^{(j)}_i(\Ab,\bb)=0$, where $j$ ranges from $1$ to $p$ and $i$ is simply an
identifying index for the individual conditions, each of which is a polynomial of
degree $j$.  For details concerning RK order
conditions, see for instance \cite{Butcher2008}.
We ensure that those conditions are satisfied up to order $p$ and we seek
to minimize the Euclidean norm of the truncation error coefficients of
order $p+1$ \cite{Butcher2008,Hairer1991}:
$$C^{(p+1)} = \left(\sum_i \tau^{(p+1)}_i \right)^{1/2}.$$

Memory requirements for RK methods are typically on the order of
$s \times \ndof$. To avoid the need for large amounts of memory, we employ the low-storage
algorithm presented in \cite{ketcheson2010a}, which reduces the requirement
to $3\times \ndof$, i.e. 3 registers per stage. The coefficients of the methods 
are provided in terms of the low-storage formulation, which is given in
Algorithm \ref{alg-3s}.
This algorithm also retains the previous solution values so that a step can
be restarted if a prescribed stability or accuracy condition is not met.

\begin{algorithm}[ht] \topcaption{Low storage implementation (3S*)} \label{alg-3s} 
\label{alg:bisection}
\begin{algorithmic}
    \State $S_3 \gets u^n$
    \State $S_2 \gets 0$
    \State $S_1 \gets u^n$ 
    \For{$i=1:s$}
        \State $t \gets t^n + c_i\ \dt$
        \State $S_2 \gets S_2 + \delta_{i}\ S_1$
        \State $S_1 \gets \gamma_{1,i}\ S_1 + \gamma_{2,i}\ S_2 + \gamma_{3,i}\ S_3 + \beta_{i}\ \dt \ F(S_1)$ 
    \EndFor
\State $u^{n+1} \gets S_1$
\end{algorithmic}
\end{algorithm}

The optimization problem may be stated formally as follows
\begin{problem}[RK method optimization] \label{prob:rkcoeff}
    \begin{equation*}
      \begin{aligned}
        \textup{Choose } \Ab,\bb \textup{ to } & \minimize && C^{(p+1)} \\
        & \textnormal{subject to} && \\
        & && \tau^{(j)}_i(\Ab,\bb) = 0 & (0 \le j \le p) \\
        & &&  \bb^T \Ab^{j-1}\eb=\beta_j & (0\le j\le s) \\
        & && \Gamma(\Ab,\bb) = 0
      \end{aligned}
    \end{equation*}
\end{problem}
Here $\Gamma(\Ab,\bb)=0$
represents the conditions necessary for the method to be written in low-storage
form. In practice, we impose those conditions implicitly by taking the
low-storage coefficients as decision variables and computing the Butcher 
coefficients $(\Ab,\bb)$ from them.

We use the RK-opt ("Runge--Kutta optimization") package to search for optimized 
methods. This software uses MATLAB's {\tt fmincon} function with 
the interior-point algorithm and the multi-start global optimization toolbox. 
Six hundred random initial guesses were used to find each optimized RK 
method. The RK-opt package and its extensions \cite{Ketcheson2012b} are freely available at 
\url{https://github.com/ketch/RK-opt}.

\subsection{Efficiency and CFL number}\label{subsec:efficiency}
Time integration with an explicit method always incurs a step size
restriction $\nu\le\dtm$ based on stability.  Accuracy typically also
leads to a constraint on the time step, which for hyperbolic problems
translates to a constraint on the CFL number, of the form $\nu\le \dta$, where $\dta$ is the
largest CFL number satisfying a prescribed error tolerance.
Other concerns, such as positivity, may further restrict the CFL number,
but we focus on $\dtm$ and $\dta$.

If $\dtm<\dta$, stability is the more restrictive concern and
the relative efficiency of two RK methods of order $p$ can be 
measured as the ratio of the maximum effective stable CFL number $\dtm/s$:
\begin{equation}\label{eq:stab-efficiency-def}
\chistab = \frac{\sigma \,  \nu_1/s_1}{\sigma \, \nu_2/s_2} = \frac{\nu_1/s_1}{\nu_2/s_2},
\end{equation}
where $\sigma$ denotes a safety factor applied to both schemes.
If $\chistab> 1$ then method 1 is more efficient. This
quantity measures the relative time interval integrated per unit work \cite{Kennedy2000a}.

On the other hand, if $\dta<\dtm$, accuracy is the more restrictive
concern so relative efficiency should be based on the ratio of step sizes,
giving an equivalent global error.  A first estimate, assuming that
local errors simply accumulate, yields the relative efficiency measure
\begin{equation}\label{eq:acc-efficiency-global-def}
\chiacc = \left(\frac{C^{(p+1)}_2}{C^{(p+1)}_1}\right)^{\frac{1}{p}}\frac{s_2}{s_1},
\end{equation}
where $C^{(p+1)}_1$ and $C^{(p+1)}_2$ are the principal error norms of 
the two RK schemes (see Section \ref{subsec:rk-coeffs}). 
Note that this measure is meaningful only if both schemes have order $p$.
Although \eqref{eq:acc-efficiency-global-def} is probably 
too simplistic because the error at each time step feeds back into the 
computation at the next step, it is used as a guideline for the selection 
of RK schemes among the optimized methods presented in this paper\footnote{Notice that \eqref{eq:acc-efficiency-global-def} differs from Equation $(28)$ 
in \cite{Kennedy2000a} by the exponent.  In Kennedy et al. 
\cite{Kennedy2000a} the power of $C^{(p+1)}_2/C^{(p+1)}_1$ is set to
$1/(p+1)$ because the control of the local truncation error 
is the main objective of the comparison \cite{Hosea1994}.}.

\begin{figure}[htbp!]
\centering
\subfigure[Optimal 2$^\mathrm{nd}$-order methods vs. ERK(2,2).]{
\includegraphics[scale=0.38]{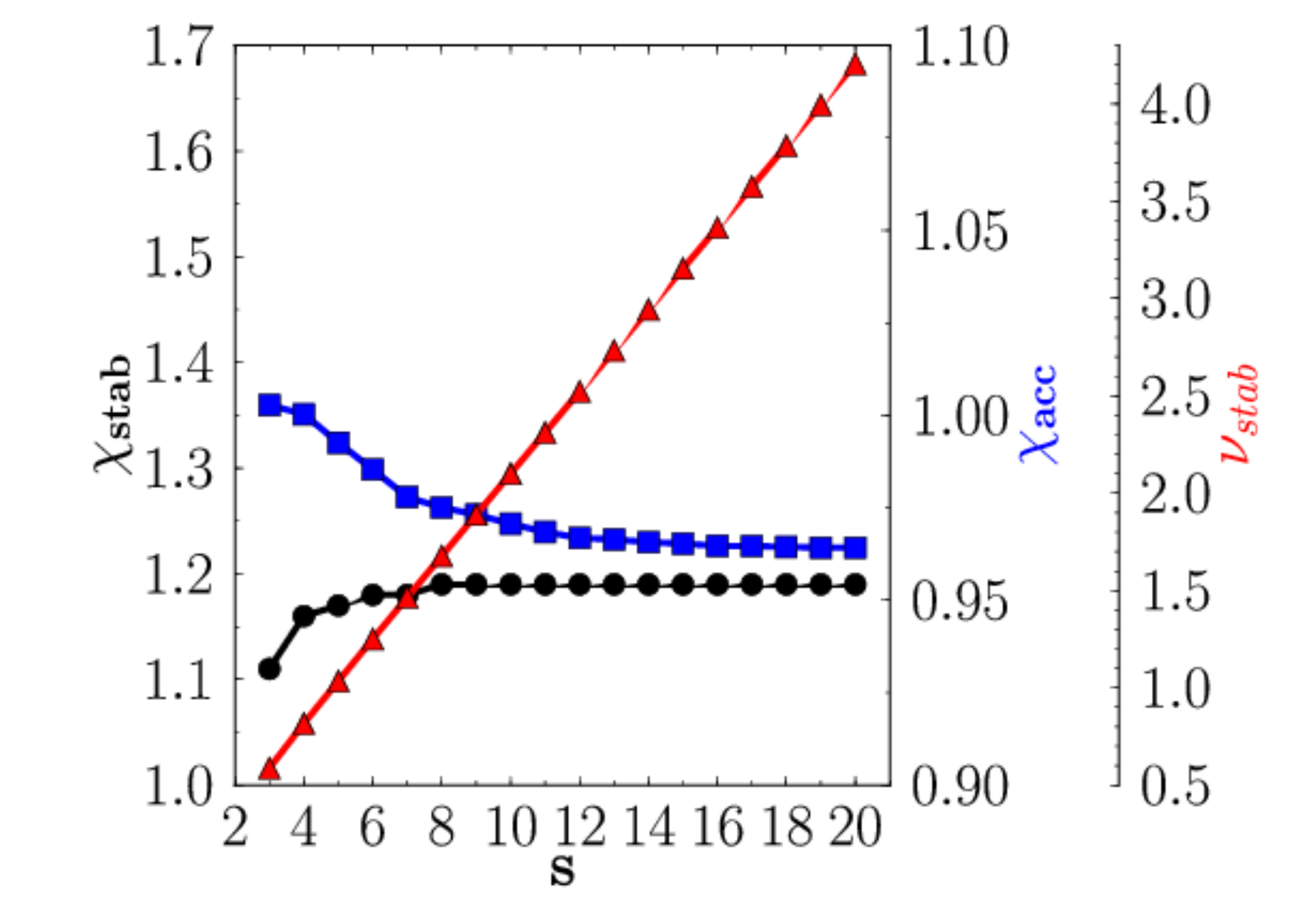}
\label{fig:efficiencies-optim-erk2}}
\hspace{-0.75cm}
\subfigure[Optimal 3$^\mathrm{rd}$-order methods vs. ERK(3,3).]{
\includegraphics[scale=0.38]{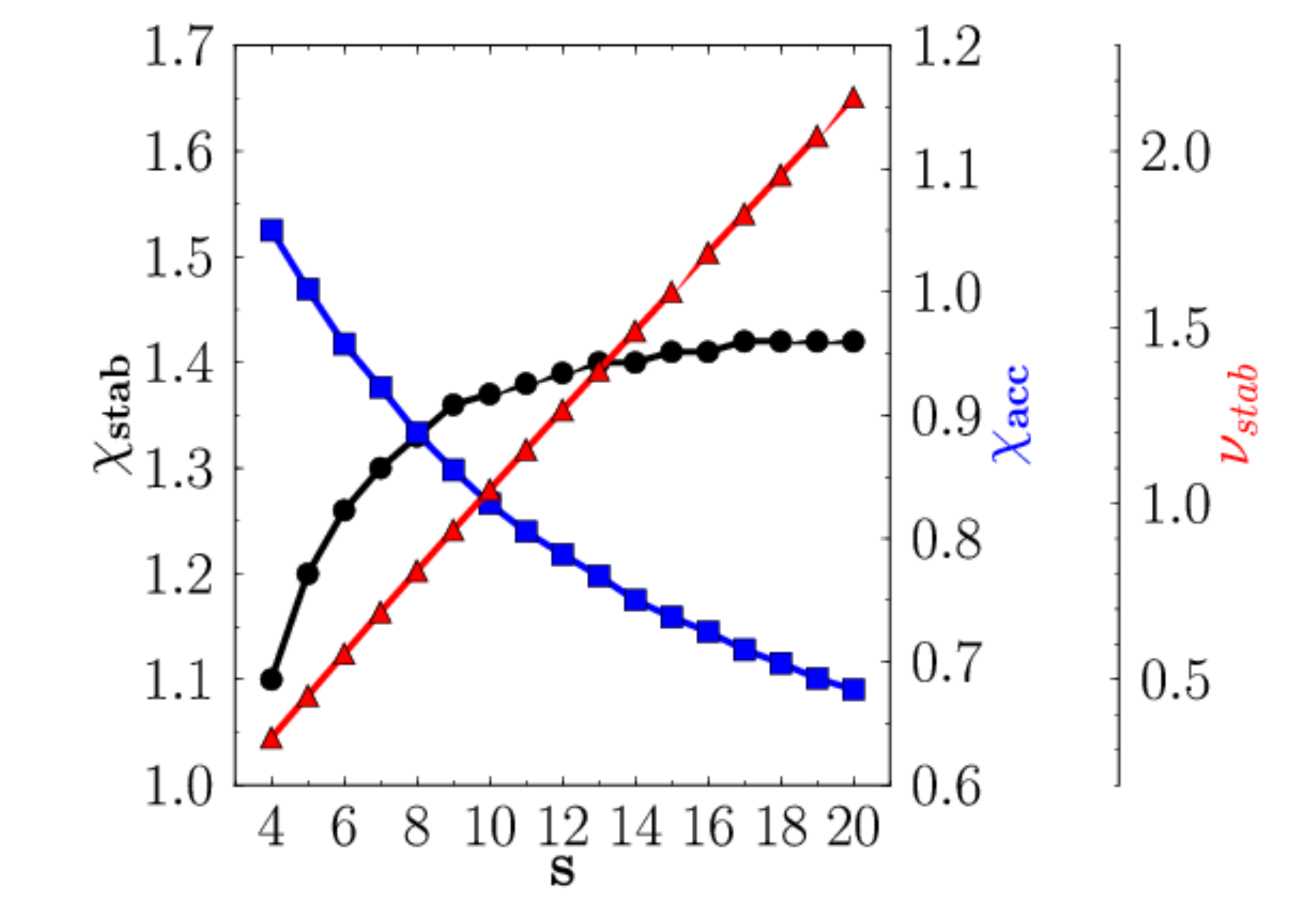}
\label{fig:efficiencies-optim-erk3}}

\vspace{1cm}

\subfigure[Optimal 4$^\mathrm{th}$-order methods vs. ERK(4,4).]{
\includegraphics[scale=0.38]{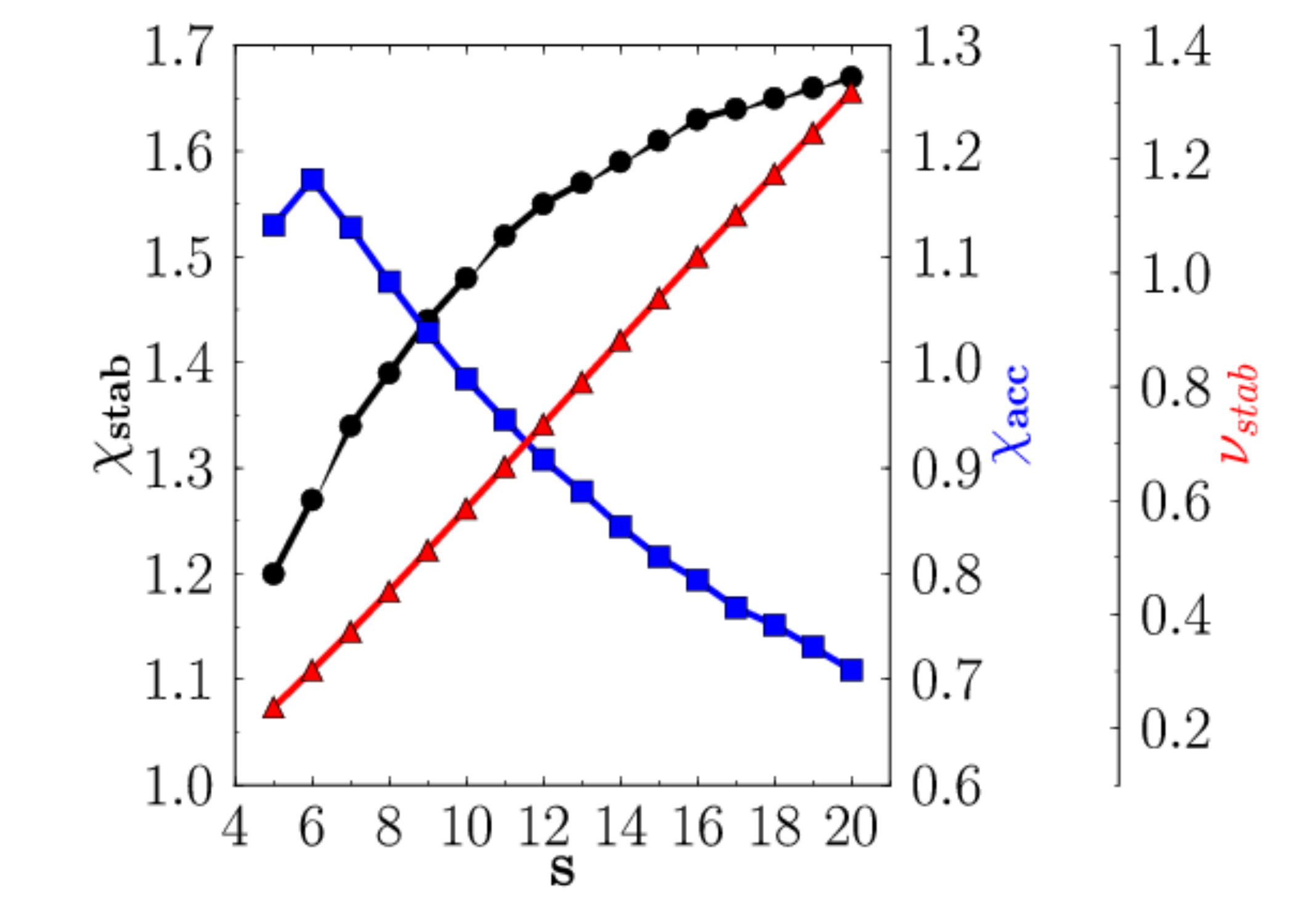}
\label{fig:efficiencies-optim-erk4}}
\hspace{-0.75cm}
\subfigure[Optimal 5$^\mathrm{th}$-order methods vs. ERKF(6,5).]{
\includegraphics[scale=0.38]{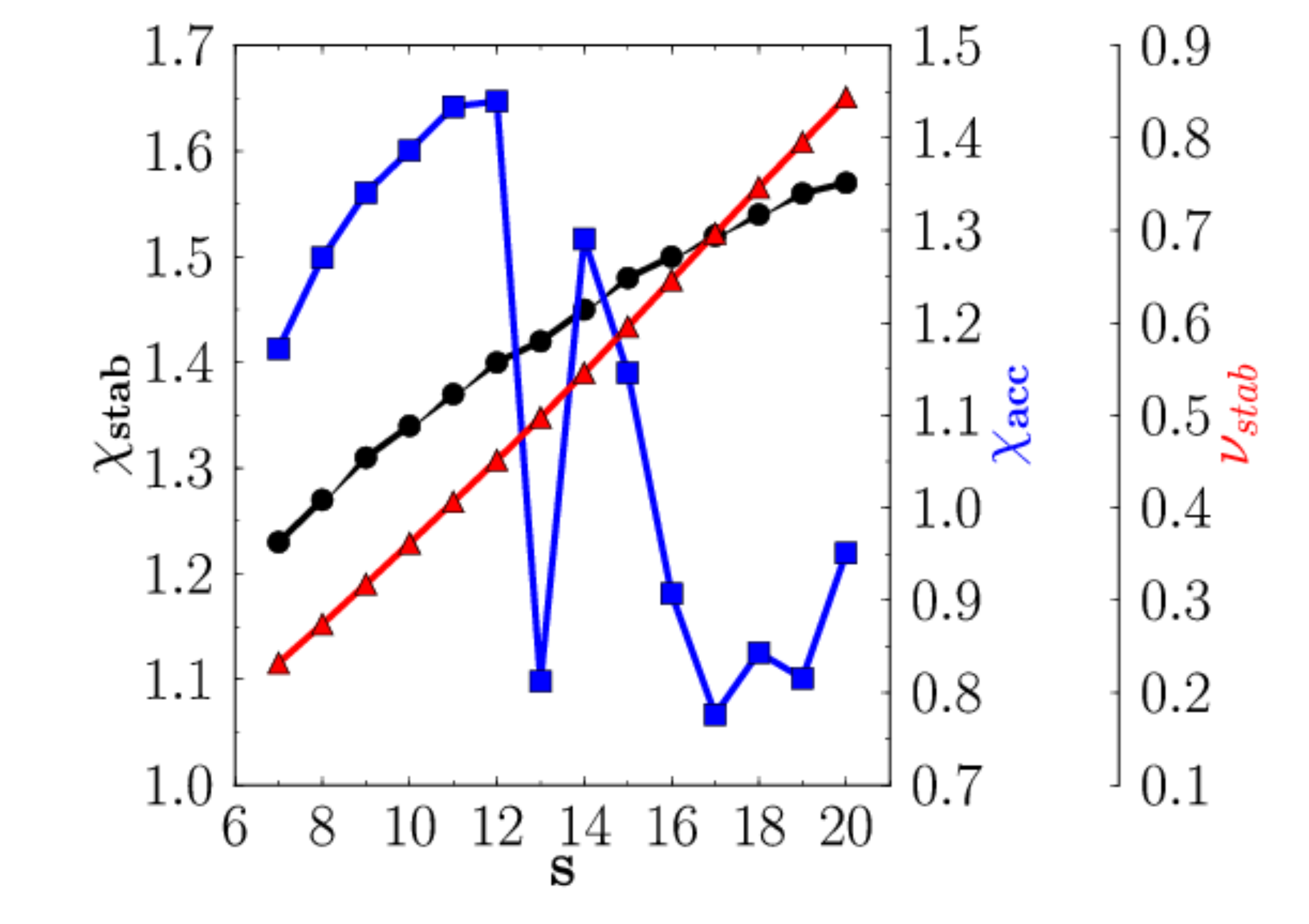}
\label{fig:efficiencies-optim-erk5}}

\caption{Efficiencies and maximum linearly stable CFL number of the optimal 
ERK methods over some traditional ERK schemes of the same 
accuracy.\label{fig:efficiencies-optim-erk}}
\end{figure}

Figure~\ref{fig:efficiencies-optim-erk} shows both $\chistab$  
and $\chiacc$ for 2$^\mathrm{nd}$- to 5$^\mathrm{th}$-order optimized schemes over widely 
used traditional explicit 
ERK methods of the same accuracy: the mid-point rule ERK(2,2); Heun's 3-stage 
3$^\mathrm{rd}$-order ERK(3,3) method \cite{Heun1900}; the classical 4-stage 
4$^\mathrm{th}$-order ERK(4,4) \cite{Ku1901}; and the 6-stage 5$^\mathrm{th}$-order Runge--Kutta--Fehlberg 
ERKF(6,5) method  \cite{Fehlberg1969}. The maximum stable CFL number $\dtm$ defined in Section 
\ref{subsub:spectrum} is also shown for completeness. 

\subsection{Discussion}\label{subsec:discussion-optimized-rk}
We provide optimized methods for $s\le 20$, because for larger values of $s$,
the convex solvers used in the algorithm of \cite{Ketcheson2012a} often fail due
to poor numerical conditioning.  
However, Figures~\ref{fig:efficiencies-optim-erk2} and \ref{fig:efficiencies-optim-erk3}
already show that the marginal efficiency gain achieved by adding another
stage, becomes vanishingly small for large $s$.

Indeed, the asymptotic efficiency gain
that could be achieved by using additional stages is bounded, since the
classical CFL theorem implies that the scheme cannot be stable for a CFL number
greater than $s$ \cite{Sanz-Serna1986a}.  
An even tighter bound can be inferred by recalling that
the stability region of an $s$-stage ERK cannot contain the closed disk with
diameter $[-2\,s,0]$ as a proper subset \cite{jeltsch1978largest}.  
By determining the largest of such disks
contained in the spectrum of the SD method, upper bounds on the efficiency
of optimized methods can be obtained. A further refinement can be 
obtained by using Theorem 5 of \cite{Sanz-Serna1986a}, which refers to ellipses
instead of only disks.

These considerations imply that, for hyperbolic PDE discretizations, 
only a limited number of stages are necessary to realize most of the
potential efficiency gain.  Based on Figure \ref{fig:efficiencies-optim-erk},
it seems that the number of stages that provides a significant improvement
increases with the order $p$.

The blue lines indicate that the global error efficiency $\chiacc$ 
of our schemes generally decreases with the number of stages.
However, for all schemes, $\chiacc$ is within $40\%$ of the reference scheme value.
We emphasize that accuracy is not the primary concern in the
design of these schemes; certainly better accuracy could be obtained if one
were willing to sacrifice some stability.

Tables \ref{tab:efficiencies-selected-optimal-erk2} to
\ref{tab:efficiencies-selected-optimal-erk5} 
list the value of $\chistab$ and $\chiacc$ of the optimized
schemes that are used for the test problems in the next section.
Two methods have been
selected for each order of accuracy. Those with fewer number of stages 
have an accuracy efficiency close to that of the reference methods, whereas the
schemes with a large number of stages are characterized by a large 
value of the stability efficiency and an accuracy efficiency which is greater 
than $0.7$. The coefficients of the selected optimized 
methods in terms of the low-storage formulation (see Algorithm \ref{alg-3s}
) are listed in Appendix \ref{app:rk_coeffs}.

\begin{table}
\centering
\begin{tabular}{c|c|cc|cc}  
Method & {$s$} & {$\dtm/s$} & $C^{(p+1)}$ & $\chistab$ & $\chiacc$ \\ \hline
Midpoint ERK(2,2) & {\em 2}& {\em 1.7678\e{-01}} & {\em 1.7180\e{-01}} & {\em 1} & \em 1 \\ \hline 
Optimal ERK(3,2) & {3} &{1.9587\e{-01}} & 7.5938\e{-02} & 1.11 & 1.00 \\ 
Optimal ERK(8,2) & {8} &{2.0968\e{-01}} & 1.1294\e{-02} & 1.19 & 0.98\\ 
\end{tabular} 
\caption{Step size, stability efficiency $\chistab$ and estimation of
accuracy efficiency $\chiacc$ of the selected optimal
2$^\mathrm{nd}$-order ERK methods.  The reference 2-stage method (corresponding to the
values in italics) is the midpoint ERK method.
\label{tab:efficiencies-selected-optimal-erk2}}
\end{table} 

\begin{table}
\centering
\begin{tabular}{c|c|cc|cc}  
Method & {$s$} & {$\dtm/s$} & $C^{(p+1)}$ & $\chistab$ & $\chiacc$ \\  \hline
Heun's ERK(3,3) & {\em 3}& {\em 7.5739\e{-02}} & \em 4.6296\e{-02}& \em 1 & \em 1 \\ \hline  
Optimal ERK(5,3) & {5} &{9.0719\e{-02}} & 9.9290\e{-03} & 1.20 & 1.00 \\  
Optimal ERK(18,3) & {17} &{1.0718\e{-01}} & 7.1115\e{-04} & 1.42 & 0.71 \\ 
\end{tabular} 
\caption{Step size, stability efficiency $\chistab$ and estimation of global accuracy efficiency $\chiacc$ of the selected optimal 3$^\mathrm{rd}$-order ERK methods. 
The reference method (corresponding to the values in italics) is Heun's 3-stage method.
\label{tab:efficiencies-selected-optimal-erk3}}
\end{table} 

\begin{table}
\centering
\begin{tabular}{c|c|cc|cc}
Method & {$s$} & {$\dtm/s$} & $C^{(p+1)}$ & $\chistab$ & $\chiacc$ \\ \hline
Kutta's ERK(4,4) & {\em 4}& {\em 3.9534\e{-02}} & {\em 1.4505\e{-02}} & {\em 1} & \em 1 \\ \hline 
Optimal ERK(9,4) & {9} &{5.6977\e{-02}} & 5.0640\e{-04} & 1.44 & 1.03 \\ 
Optimal ERK(18,4) & {18} &{6.5233\e{-02}} & 1.1087\e{-04} & 1.65 & 0.75\\ 
\end{tabular} 
\caption{Step size, stability efficiency $\chistab$ and estimation of
accuracy efficiency $\chiacc$ of the selected optimal
4$^\mathrm{th}$-order ERK methods.  The reference 4-stage method (corresponding to the
values in italics) is the classic ERK method.
\label{tab:efficiencies-selected-optimal-erk4}}
\end{table} 

\begin{table}
\centering
\begin{tabular}{c|c|cc|cc}
Method & {$s$}  & {$\dtm/s$} & $C^{(p+1)}$ & $\chistab$ & $\chiacc$ \\ \hline
Fehlberg ERK(6,5) & {\em 6}& {\em 2.6916\e{-02}} & {\em 3.3557\e{-03}} & {\em 1} & \em 1 \\ \hline 
Optimal ERK(10,5) & {10}   &{3.6164\e{-02}} & 5.0975\e{-05} & 1.34 & 1.39 \\ 
Optimal ERK(20,5) & {20}   &{4.2195\e{-02}} & 1.0490\e{-05} & 1.57 & 0.95\\ 
\end{tabular} 
\caption{Step size, stability efficiency $\chistab$ and estimation of
accuracy efficiency $\chiacc$ of the selected optimal
5$^\mathrm{th}$-order ERK methods.  The reference method (corresponding to the
values in italics) is Fehlberg's 6-stage, 5$^\mathrm{th}$-order method.
\label{tab:efficiencies-selected-optimal-erk5}}
\end{table}

The classical linear stability analysis describes the growth of
truncation errors from one step to the next, but ignores the effect
that roundoff and truncation errors in intermediate stages may have
within a single step.
Although, for conventional ERK methods, the accumulation of 
round-off errors during a single time step is negligible, it must be taken into
account for schemes with a large number of stages. In fact, in the application 
of ERK schemes with many stages to time dependent PDEs, there can be a serious 
accumulation of errors that may even render methods unusable; this is referred
to as internal instability \cite{verwer1996explicit,Verwer1990}.
Since most of our new schemes 
use many stages, a thorough analysis of their internal stability  
properties has also been performed. All the schemes are internally stable.

\section{Applications}\label{sec:applications}
In order to asses the efficiency and the accuracy of our new 
ERK schemes, we have performed a series of numerical simulations. 
The computations run on a machine with $2 \times 2.4 $ GHz Quad-Core 
Intel Xeon, using the Coolfluid 3 collaborative simulation
environment \cite{cf3}. 
Sixteen gigabytes of RAM were available. The grids have been generated using
Gmsh software \cite{Geuzaine2009}.

In Coolfluid 3, the CFL-number in two dimensions is defined as
\begin{equation}
  \nu = \dt\ \left( \frac{u}{\Delta x} + \frac{v}{\Delta y}  \right)\ ,
\end{equation}
with $u$ and $v$ the x- and y-components of the wave speed $\vec{a}$, and with $\Delta x$ and $\Delta y$ the width and height of a Cartesian grid cell. This definition, when applied to the Roe scheme on a structured Cartesian mesh leads to
a positive and stable discretization for a CFL number smaller than unity \cite{Deconinck1996}. 
In practical computations it has been observed that slightly larger time steps may
be used without affecting stability \cite{Deconinck1996}.

\subsection{Order verification}
\begin{figure}[htbp!]
\centering
\subfigure[2$^\mathrm{nd}$-order ERK methods.]{
\includegraphics[scale=0.38]{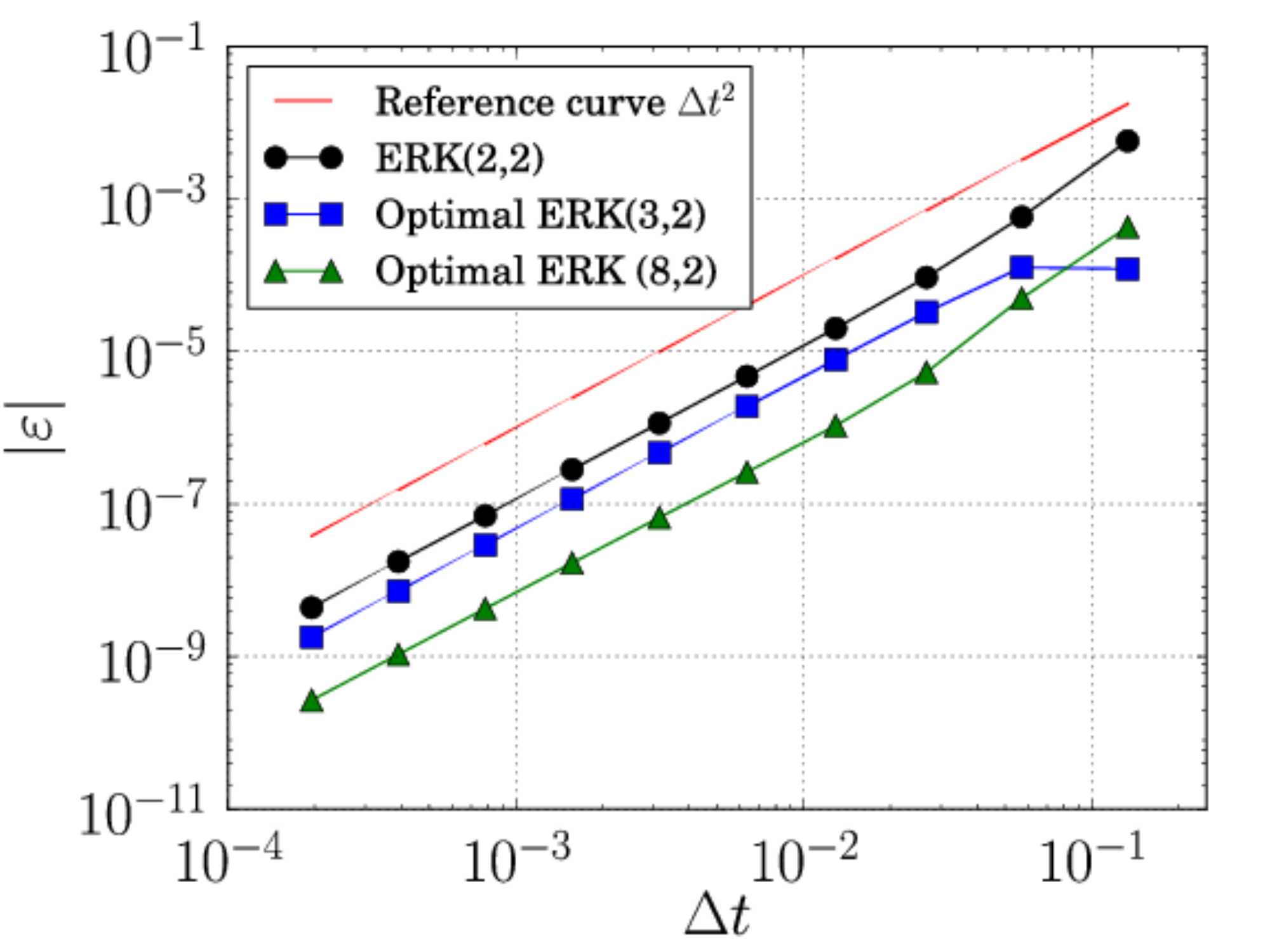}
\label{fig:convergence-optim-erk2}}
\hspace{-0.55cm}
\subfigure[3$^\mathrm{rd}$-order ERK methods.]{
\includegraphics[scale=0.38]{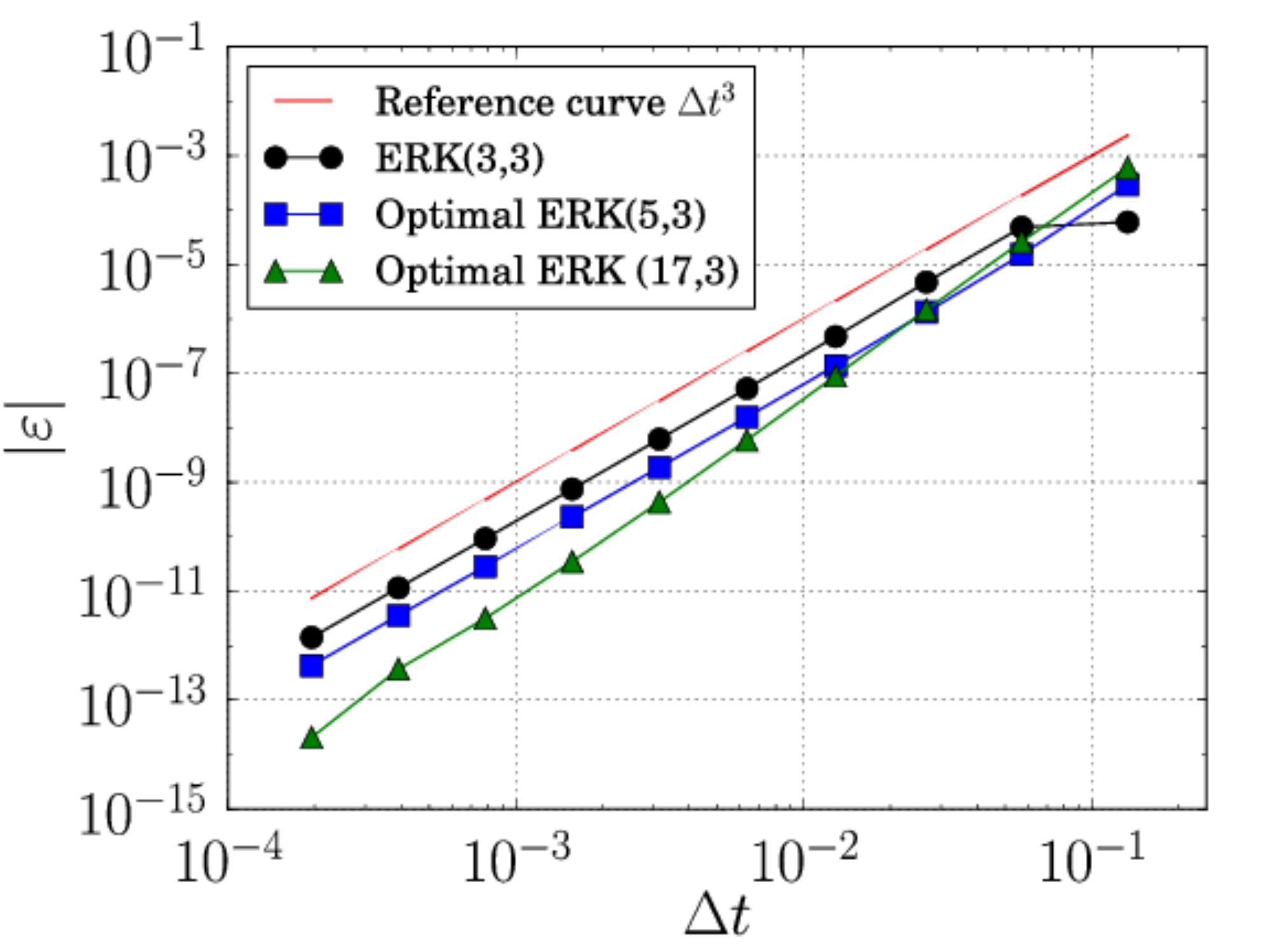}
\label{fig:convergence-optim-erk3}}

\vspace{1cm}

\subfigure[4$^\mathrm{th}$-order ERK methods.]{
\includegraphics[scale=0.38]{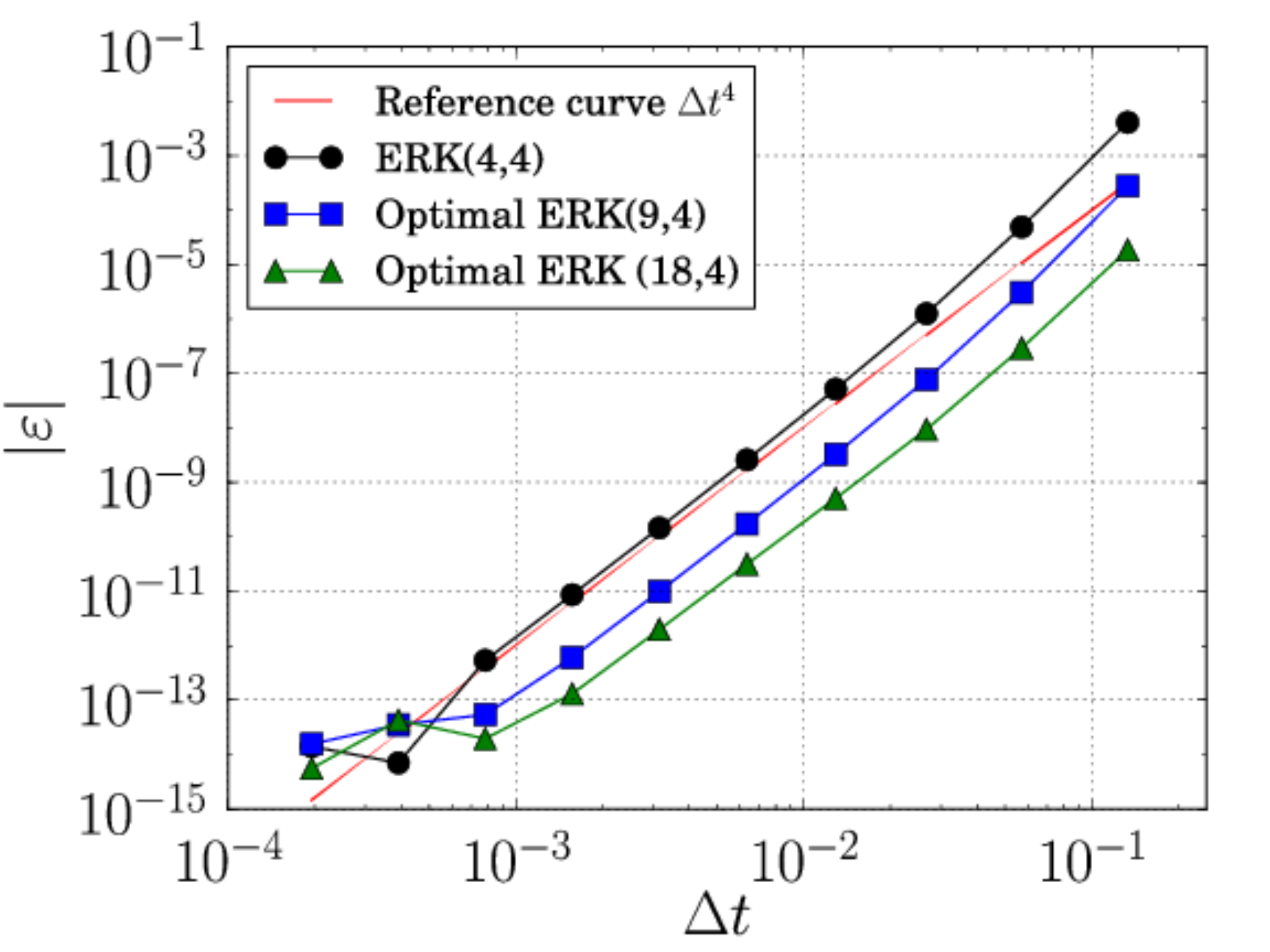}
\label{fig:convergence-optim-erk4}}
\hspace{-0.55cm}
\subfigure[5$^\mathrm{th}$-order ERK methods.]{
\includegraphics[scale=0.38]{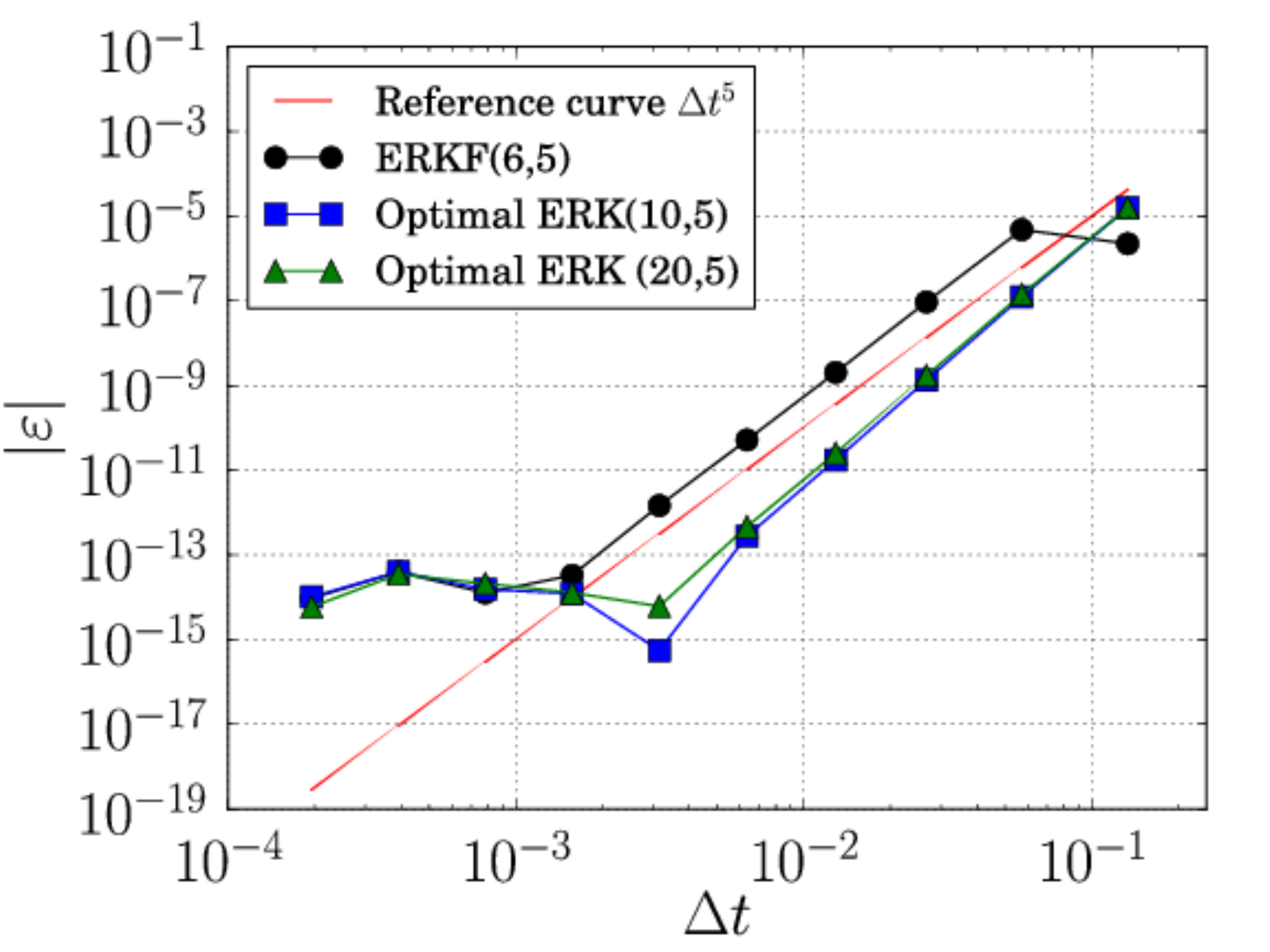}
\label{fig:convergence-optim-erk5}}

\caption{Convergence study of the optimized ERK methods.\label{fig:convergence}}
\end{figure}

In this section we present the convergence study of the optimized ERK 
scheme listed in Tables \ref{tab:efficiencies-selected-optimal-erk2} to
\ref{tab:efficiencies-selected-optimal-erk5}. We integrate a system of nonlinear
non-autonomous system of first order ODE 
\cite{Stanescu1998,Niegemann2011}
\begin{equation} 
    \begin{split}
\frac{d q_1}{d t} &= \frac{1}{q_1} - q_2\,\frac{e^{t^2}}{t^2} - t  \\ 
\frac{d q_2}{d t} &= \frac{1}{q_2} - e^{t^2} - 2\,t\,e^{-t^2} , 
    \end{split}
\label{eq:stanescu1998}
\end{equation}
with the time $t$ ranging from $t^0=1$ to $t^{e}=1.4$, and with the following initial 
conditions: $q_1(t^0) = 1, \, q_2(t^0) = e^{-1}$.
The analytical solution of this system is $q_1(t) = 1/t, \,q_2(t) = e^{-t^2}$.
We use the norm of the error
\begin{equation*}
|\varepsilon(t^{e})| = |\left(Q_1(t^{e})-q_1(t^{e})\right) + \left(Q_2(t^{e})-q_2(t^{e})\right)|
\end{equation*}
to study the time integration error. Here $Q_1$ and $Q_2$ denote the numerical 
solutions. Figure \ref{fig:convergence} shows the norm of the error $|\varepsilon|$ as 
a function of the time step $\dt$. It can be 
seen that for all ERK schemes the expected order of accuracy is achieved. 
Moreover, the new optimized
methods show significantly smaller errors than the reference methods,
as expected based on their smaller error constants.

\subsection{Advection of a Gaussian wave in an annulus}
The second problem we consider is the advection of a Gaussian 
wave in a 2D annulus. Such a problem 
models the transport of a scalar conserved variable $q$ with variable advection 
speed $\vec{a}$. The conserved variables and the convective flux are then
\begin{equation}
\begin{aligned}
  \mathbf{q} &= q, \\
  \vec{\mathbf{f}} & = q.
\label{eq: linadvannulus}
\end{aligned}
\end{equation}
Therefore, the conservation law reads
\begin{equation}
	\frac{\partial q}{\partial t} + \vec{a} \cdot \vec{\nabla} q = 0,
\label{eq:lin-adv-annulus}
\end{equation}
where in a 2D Cartesian space
\begin{equation*}
\vec{a}\left(x,y\right) = \left[\begin{array}{c}u \\ v\end{array}\right] = \omega\left[\begin{array}{c}-y \\ x\end{array}\right].
\end{equation*}
Here $\omega$ is the angular velocity which is set to $\omega = 2\pi$. Note that
the ERK schemes have been optimized for the 2D advection equation with 
constant convective velocity (see \eqref{eq: linadv}) whereas here a variable velocity is used. Therefore,
the numerical results and performance presented in this section can already be used to 
partially asses the robustness of the new time stepping methods. All schemes
were also tested using a uniform advection velocity on a uniform cartesian grid.
The results of those tests are omitted since they are very similar to those
of the more challenging test problem we now consider.

The initial Gaussian wave condition  is centered at 
$x_c = 0.0$, $y_c = 7.5$ and is defined as
\begin{equation}\label{eq:GW_annulus}
q^0(x,y) = e^{-\frac{(x-x_c)^2+(y-y_c)^2}{2 \,b^2}},
\end{equation}
where the radius is set to $b = 0.6$.

\begin{figure}[ht!]
\centering
\includegraphics[scale=0.5]{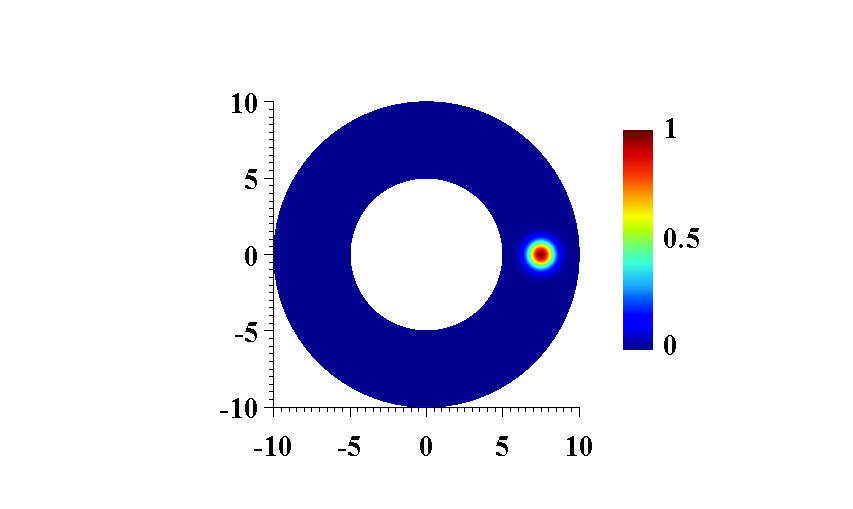}
\caption{Gaussian wave advected in the annulus at $t = 0.25$; solution computed 
with the 4$^\mathrm{th}$-order SD method and the optimal ERK(18,4) scheme.}
\label{fig:annulus_contour_4th_order}
\end{figure}

The annulus
is characterized by an internal radius $r_i = 5$ and an external radius $r_o = 10$.
One-quarter of the annulus is discretized for the actual computations.
Simulations are performed using 2$^\mathrm{nd}$- to 5$^\mathrm{th}$-order spatial and temporal discretizations
from $t^0 = 0$ to $t^{e}=0.25$ (see Figure \ref{fig:annulus_contour_4th_order}). 
Several CFL numbers ranging from $0.05$ to the maximum linearly stable one are used. 
A mesh with $110 \times 110$ (radial direction $\times$ azimuthal direction) quadrilateral
cells with a maximum aspect ratio of $1.8$ are used for the second-order 
computations. Such a mesh leads to a total number of DOFs is $48400$ which
is held constant for higher order accurate simulations by coarsening the grid
in both directions. Extrapolation boundary conditions are imposed on both
circular boundaries.

The exact solution is just a $90^\circ{}$ clock-wise roto-translation of the initial 
solution and is given by (\ref{eq:GW_annulus}) with $x_c = 7.5$, 
$y_c = 0.0$.

\begin{figure}[htbp!]
\centering
\subfigure[Error versus CPU time.]{
\includegraphics[width=0.49\linewidth]{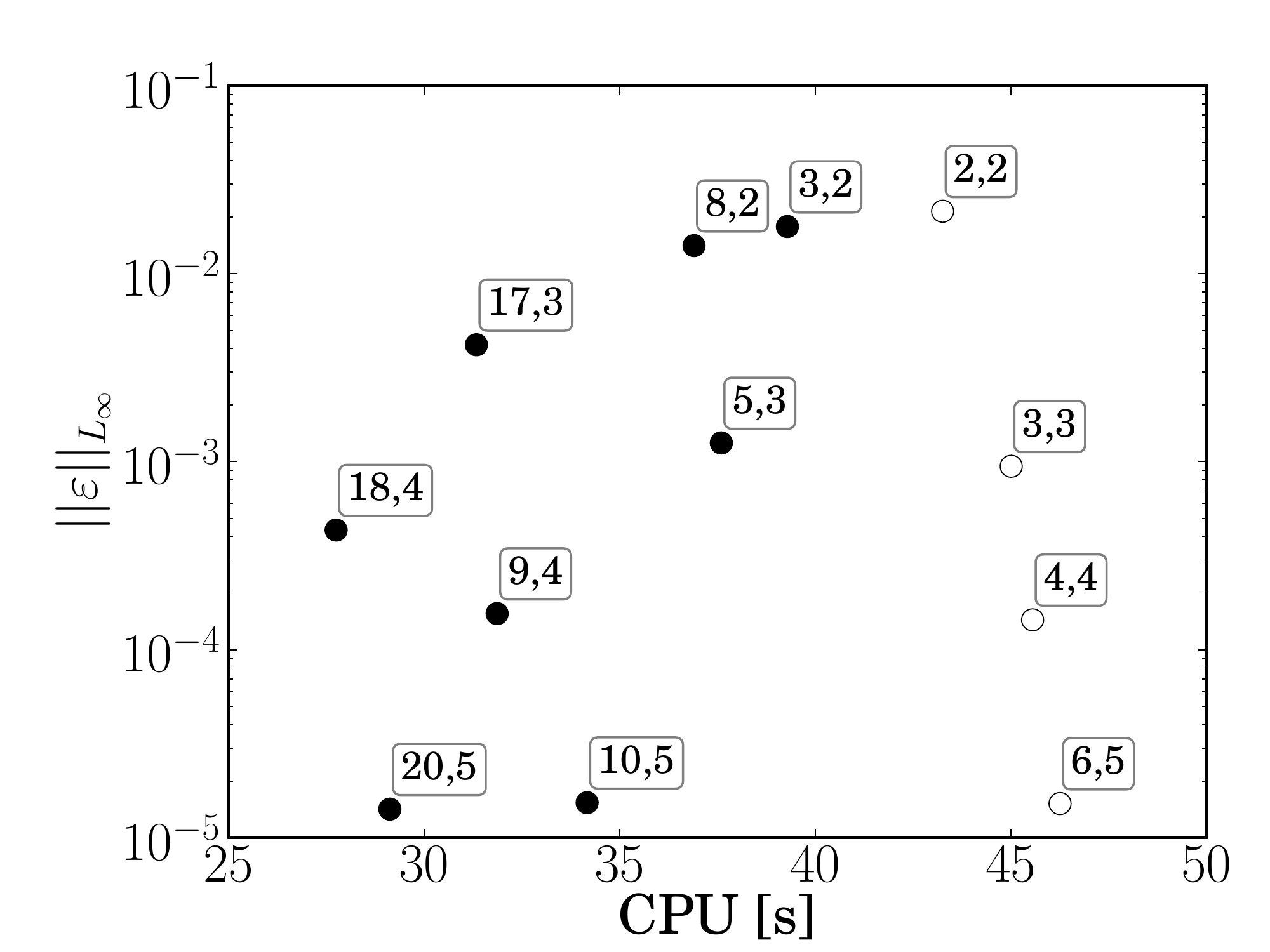}
\label{fig:error-cpu-annulus}}
\hspace{-0.72cm}
\subfigure[CPU time for each simulation.]{
\includegraphics[width=0.49\linewidth]{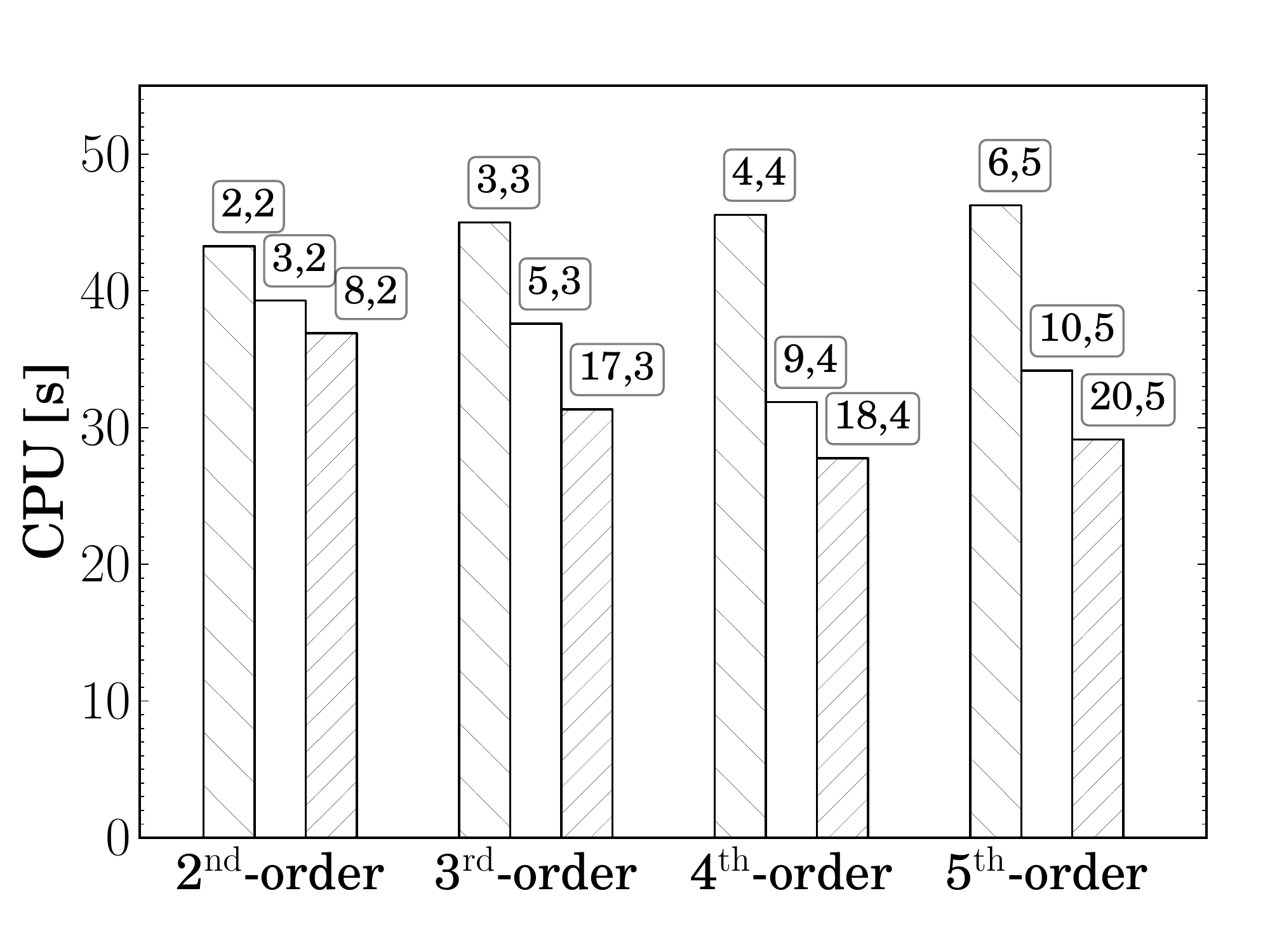}
\label{fig:cpu-annulus}}
\caption{Error and CPU time for the advection problem.
The label $s,p$ for each point indicate the number of stages $s$ and the order $p$
of the corresponding scheme.  Open circles are used for the reference methods;
closed circles are used for the optimized methods.\label{fig:err-cpu-sub-annulus}}
\end{figure}

We solve this problem with each of the reference and optimized ERK schemes 
using the predicted maximum stable CFL number for each scheme.  
In every case, the resulting computation is stable,
confirming the theoretical prediction.
Figure \ref{fig:err-cpu-sub-annulus} shows the maximum norm of the error vector 
\begin{equation}\label{eq:max_error_norm}
||\boldsymbol{\varepsilon}(t^{e})||_{L_{\infty}} = \max|\varepsilon_i| = | \max|\left(Q_i(t^{e})-q_{ex,i}(t^{e})\right)| \quad \mathrm{for} \quad i=1,2,\ldots,\ndof
\end{equation}
and the CPU time for each scheme.
Although the number of DOF is the same in all the simulations,
the error decreases rapidly with increasing order of the discretization.
Remarkably, a unit increment of the order of accuracy 
leads to a reduction of the error 
of one order of magnitude and to a faster simulation (Figure \ref{fig:error-cpu-annulus}). 
This shows the benefit of using high-order accurate methods for 
wave propagation problems. 
As predicted, some of the highly optimized ERK schemes yield somewhat larger errors.
Figure \ref{fig:cpu-annulus} highlights the speed-up obtained with 
the optimized RK schemes over the standard methods for high-order accurate 
simulations. Indeed, for 4$^\mathrm{th}$- and 5$^\mathrm{th}$-order computations
the new schemes reduce the computational time by $40\%$ and $38\%$,
respectively. These values match very well the theoretical results shown in the 
previous sections.

\begin{figure}[htbp!]
\centering
\subfigure[Optimal 2$^\mathrm{nd}$-order methods vs. ERK(2,2).]{
\includegraphics[scale=0.375]{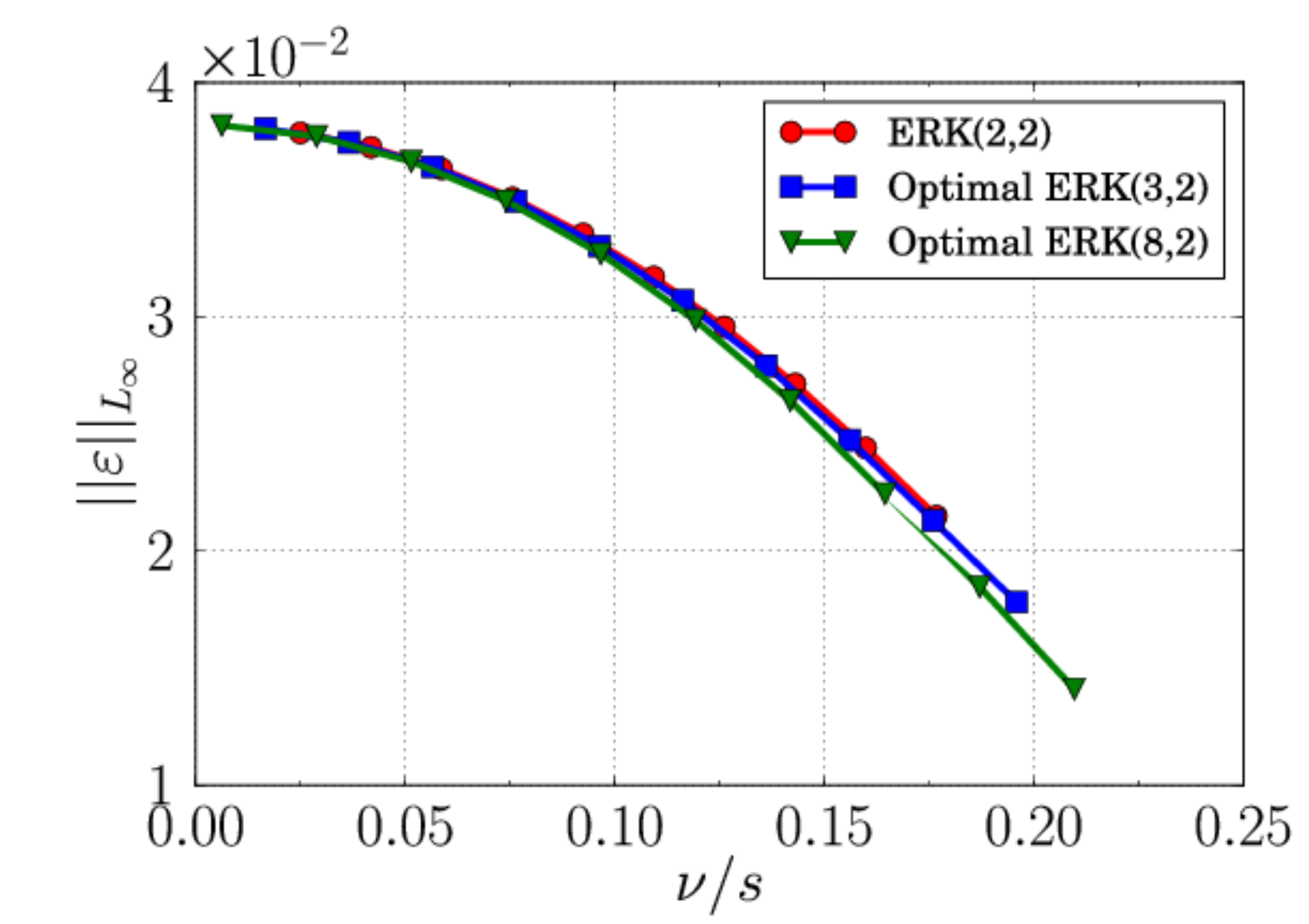}
\label{fig:annulus-optim-erk2}}
\hspace{-0.57cm}
\subfigure[Optimal 3$^\mathrm{rd}$-order methods vs. ERK(3,3).]{
\includegraphics[scale=0.375]{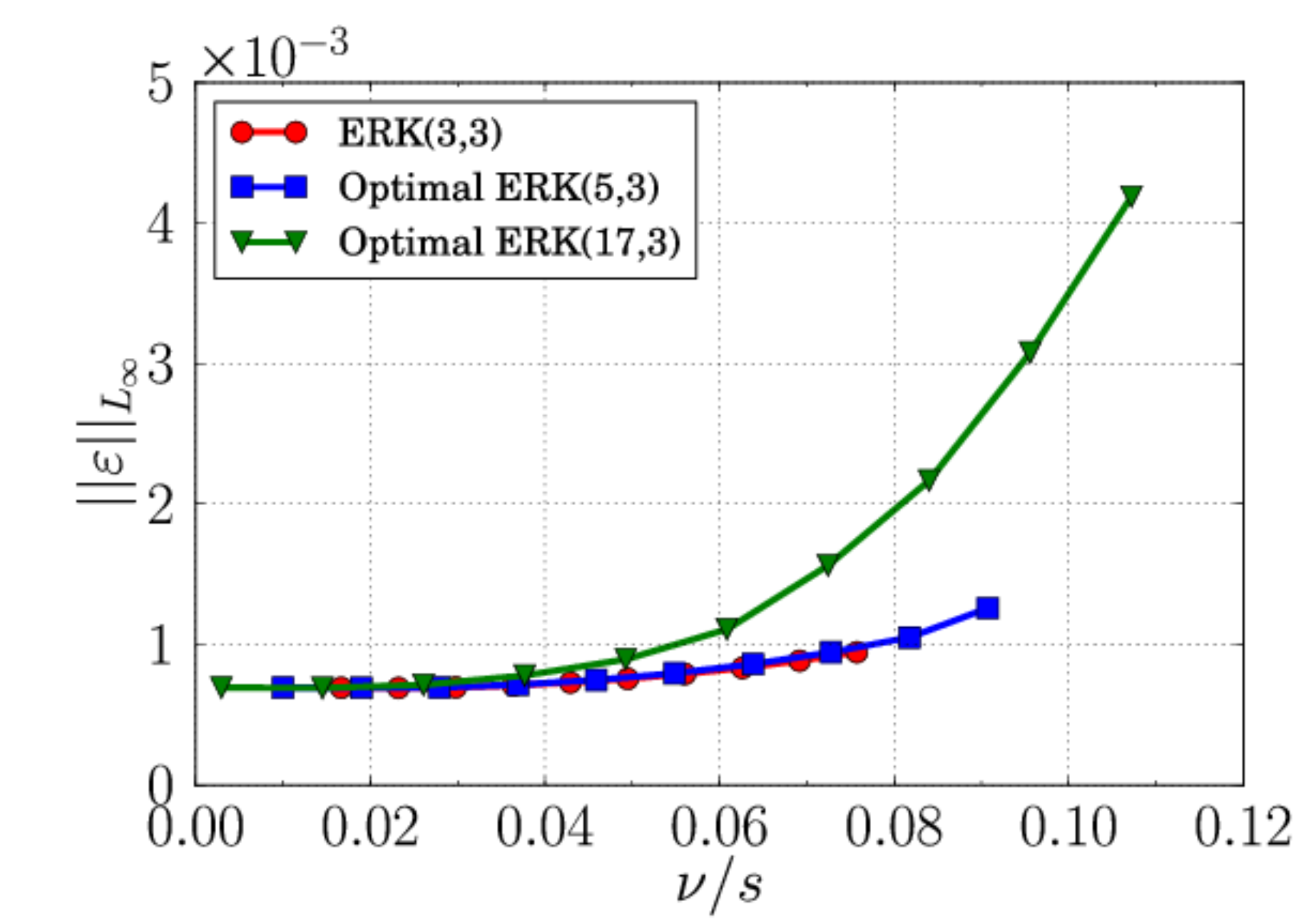}
\label{fig:annulus-optim-erk3}}

\subfigure[Optimal 4$^\mathrm{th}$-order methods vs. ERK(4,4).]{
\includegraphics[scale=0.375]{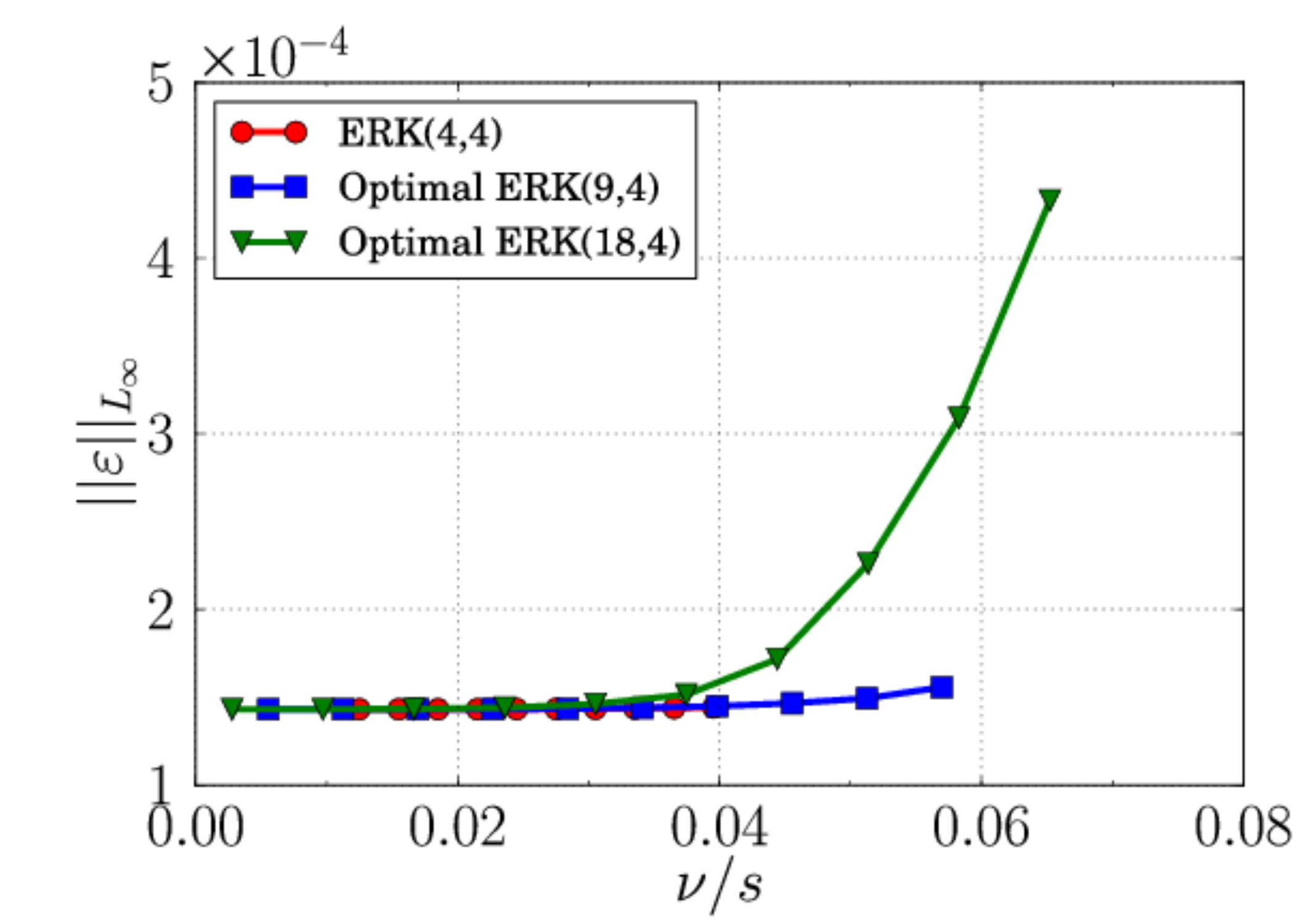}
\label{fig:annulus-optim-erk4}}
\hspace{-0.57cm}
\subfigure[Optimal 5$^\mathrm{th}$-order methods vs. ERK(6,5).]{
\includegraphics[scale=0.375]{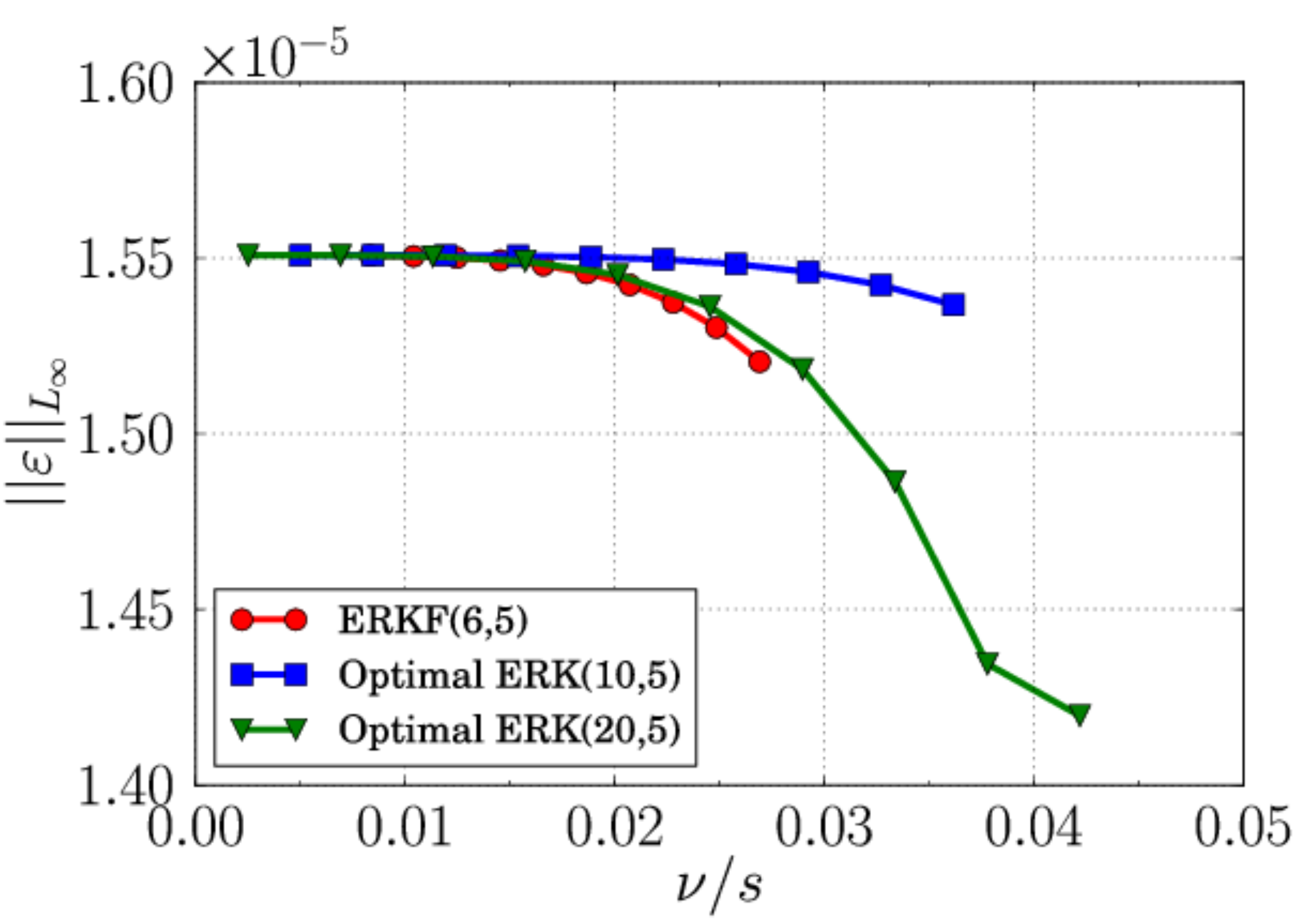}
\label{fig:annulus-optim-erk5}}

\caption{Influence of the CFL number on the maximum error norm
of a 2D Gaussian wave advecting in an annulus.\label{fig:annulus}}
\end{figure}

Figure~\ref{fig:annulus} shows the maximum norm error
as a function of the one-step effective CFL number. Here $q_{ex,i}$ 
is the exact solution at the solution point (or DOF) $i$. 

Interestingly, for some methods the error increases with increasing CFL number,
while for others it actually decreases.  The latter behavior is reminiscent
of the behavior of many low-order schemes that are more accurate for CFL numbers
close to 1 and more dissipative for small CFL numbers.
In addition,
we point out that a combination of a quasi uniform grid and the maximum 
CFL numbers obtained during the first optimization step results in stable 
full discretizations. 

\subsection{Acoustic wave propagation}
In this example we solve the linearized Euler equations (LEE), which model the propagation of small 
perturbations in a mean flow field. They are frequently used to compute 
the propagation of acoustic waves in the absence of acoustic 
sources, e.g. turbulence production. They have been successfully used to solve in
a hybrid approach for cavity flow \cite{ClausWagner2007}, jet noise \cite{Bogey2002}, 
and vortex–blade interaction \cite{Chen2005}.

The LEE are derived from the compressible Euler equations which mathematically
describe the three physical conservation laws (i.e. conservation of mass, 
conservation of momentum and conservation of energy) for an
inviscid fluid. Thus, the definitions of the conserved variables $\mathbf{q}$ and 
the flux vector $\fbv=\left[\mathbf{f}\;\mathbf{g}\;\mathbf{h}\right]^T$ 
are
\begin{subequations} \label{eq:Euler}
  \begin{align}
    \mathbf{q} & = \left(\begin{array}{c}\rho\\ \rho u\\ \rho v \\ \rho w \\ \rho E\end{array}\right), 
\label{eq:consVar} \\
 \mathbf{f} = \left(\begin{array}{c}\rho u \\ \rho u^2+p\\ \rho u v \\ \rho u w \\ u \left(\rho E+p\right)\end{array}\right), \hspace{0.1in}
 \mathbf{g} & = \left(\begin{array}{c}\rho v \\ \rho u v\\ \rho v^2+p \\ \rho v w \\ v \left(\rho E+p\right)\end{array}\right), \hspace{0.1in} \mathbf{h} = \left(\begin{array}{c}\rho w \\ \rho u w\\ \rho v w \\ \rho w^2 + p \\ w \left(\rho E+p\right)\end{array}\right).
 \label{eq:FluxComps}
  \end{align}
\end{subequations}
In these equations, $\rho$ is the mass density, $u$, $v$ and $w$ are the 
Cartesian velocity components, $p$ is thermodynamic pressure, and $E$ is 
specific total energy. The velocity vector $\vec{u}$ is $\left[u\;v\;w\right]^T$ 
and its magnitude is denoted by $\left|\vec{u}\right|$. For an ideal gas, 
which approximates well the thermodynamic
behavior of air in a wide range of thermodynamic conditions, the
specific total energy 
$E$ is related to the pressure and the velocity field by
\begin{equation}
  E= \frac{1}{\gamma-1}\frac{p}{\rho}+\frac{u^2+v^2+w^2}{2},
\label{eq:energyIdealGas}
\end{equation}
where $\gamma=1.4$ is the heat capacity ratio for air. Equation \eqref{eq:energyIdealGas} closes the hyperbolic system \eqref{eq:Euler} of five 
nonlinear PDEs with five unknowns.

The LEE are obtained from \eqref{eq:Euler} by decomposing the 
primitive flow variables $\rho$, $\vec{u}$ and $p$ into a mean flow value 
$(\cdot)_0$ and a perturbation about this mean flow $(\cdot)'$:
\begin{equation}
\begin{array}{rcl}
  \rho    &=& \rho_0    + \rho',\\
  \vec{u} &=& \vec{u}_0 + \vec{u}',\\
  p       &=& p_0       + p'.
\end{array}
\end{equation}

Substituting these relations in \eqref{eq:Euler}, subtracting the mean flow 
terms and neglecting products of perturbations, the following sets of conserved
variables $\qb$ and the flux components $\fbv=\left[\mathbf{f}\;\mathbf{g}\;\mathbf{h}\right]^T$
are obtained
\begin{equation}
  \mathbf{q} = \left(\begin{array}{c}\rho'\\ \rho_0 u'\\ \rho_0 v' \\ \rho_0 w' \\ p'\end{array}\right),\hspace{0.25in}
  \label{eq: consVarLEE}
\end{equation}

\begin{equation}
  \mathbf{f} = \left(\begin{array}{c}\rho_0 u'+u_{0}\rho' \\ \rho_0 u_{0}u'+p'\\ \rho_0 u_{0}v'    \\ \rho_0 u_{0}w'    \\ u_{0}p'+\gamma p_0 u'\end{array}\right),\hspace{0.1in}
  \mathbf{g} = \left(\begin{array}{c}\rho_0 v'+v_{0}\rho' \\ \rho_0 v_{0}u'   \\ \rho_0 v_{0}v'+p' \\ \rho_0 v_{0}w'    \\ v_{0}p'+\gamma p_0 v'\end{array}\right),\hspace{0.1in}
  \mathbf{h} = \left(\begin{array}{c}\rho_0 w'+w_{0}\rho' \\ \rho_0 w_{0}u'   \\ \rho_0 w_{0}v'    \\ \rho_0 w_{0}w'+p' \\ w_{0}p'+\gamma p_0 w'\end{array}\right).
\label{eq: lineuler2}
\end{equation} 

This procedure also leads to a source term involving mean flow gradients (right-hand side of 
\eqref{eq:genHypSys}): 
\begin{equation}
  \mathbf{s}={-\left(\begin{array}{c}0\\ \left(\rho_0 u'+u_{0}\rho'\right)\frac{\partial u_{0}}{\partial x}+\left(\rho_0 v'+v_{0}\rho'\right)\frac{\partial u_{0}}{\partial y}+\left(\rho_0 w'+w_{0}\rho'\right)\frac{\partial u_{0}}{\partial z}\\ \left(\rho_0 u'+u_{0}\rho'\right)\frac{\partial v_{0}}{\partial x}+\left(\rho_0 v'+v_{0}\rho'\right)\frac{\partial v_{0}}{\partial y}+\left(\rho_0 w'+w_{0}\rho'\right)\frac{\partial v_{0}}{\partial z}\\ \left(\rho_0 u'+u_{0}\rho'\right)\frac{\partial w_{0}}{\partial x}+\left(\rho_0 v'+v_{0}\rho'\right)\frac{\partial w_{0}}{\partial y}+\left(\rho_0 w'+w_{0}\rho'\right)\frac{\partial w_{0}}{\partial z}\\ \left(\gamma-1\right)\left(p'\left(\frac{\partial u_{0}}{\partial x}+\frac{\partial v_{0}}{\partial y}+\frac{\partial w_{0}}{\partial z}\right)-u'\frac{\partial p_0}{\partial x}-v'\frac{\partial p_0}{\partial y}-w'\frac{\partial p_0}{\partial z}\right)\end{array}\right)},
\end{equation}
which partially accounts for the refraction effects. The source term is zero in 
case of a uniform mean flow.

\begin{figure}[ht!]
\centering
\includegraphics[scale=0.5]{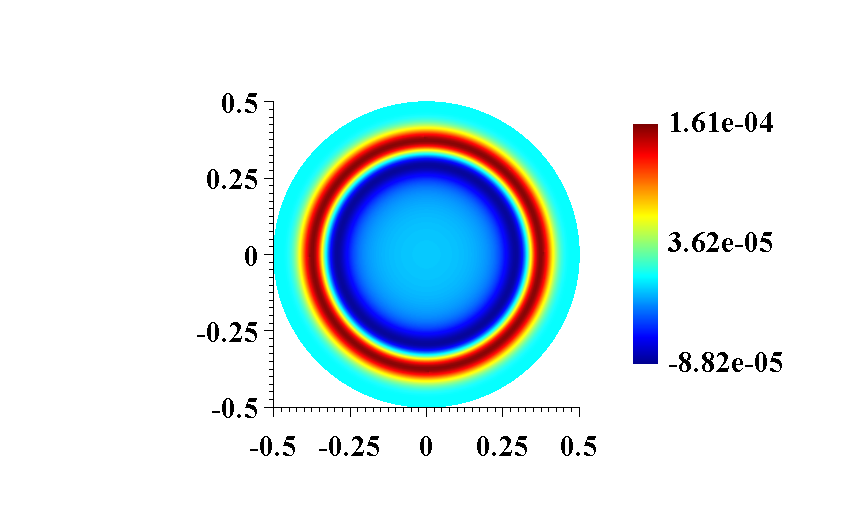}
\caption{Acoustic pressure contour $p'$ at $t = 0.3$; solution computed with the 
4$^\mathrm{th}$-order SD method and the optimal ERK(18,4) scheme.}
\label{fig:acoustic_pressure_contour_4th_order}
\end{figure}

The initial solution for this numerical test has a Gaussian profile centered at the origin of the 
axes and it is given by
\begin{equation}
\begin{array}{rcl}
  \rho'    &=& 10^{-3} e^{-\frac{x^2 + y^2}{b^2}},\\
  p'       &=& \cb^2 \ \rho',\\
  u'       &=& v' = 0,
\end{array}
\end{equation}
where the radius of the Gaussian pulse is set to $b=0.05$. The uniform mean flow 
variables are $\rhob = 1$, $\pb = 1$, $\gamma = 1.4$ and $u_0=v_0=0$. Simulation
are done from $t^0= 0$ to $t^e = 0.3$.
The numerical domain is a circle with radius $r = 0.5$, which is also 
centered at the origin of the axes. For the 2$^\mathrm{nd}$-order calculations, a
mesh with $180 \times 45$ (radial direction $\times$ azimuthal direction) quadrilateral
cells with a maximum aspect ratio of $1.85$ is used. Therefore the total number of 
DOFs is $32400$. As for the previous numerical test, this number
is kept constant for higher order accurate simulations by coarsening the grid 
in both directions. Simple extrapolation boundary conditions are used.
Figure \ref{fig:acoustic_pressure_contour_4th_order} shows the contour plot of 
the acoustic pressure field at $t = t^e = 0.3$. 

The exact solution for the acoustic pressure field $p'=p-\pb$ obtained by 
integrating the LEE is used as a reference solution to compute the numerical 
error. Its analytical expression is given by \eqref{eq:exact_acustic_pulse}
\begin{equation}\label{eq:exact_acustic_pulse}
  p'\left(t,x,y\right) = \frac{\cb^2b^2}{2} \int_0^{+\infty}e^{-\left(\frac{\xi b}{2}\right)^2}\cos\left(\xi \cb t\right)J_0\left(\xi\eta\right)\xi d\xi,
\end{equation}
with $\eta = \sqrt{\left(x-0.5\right)^2+\left(y-0.5\right)^2}$ and $J_0$ the 
Bessel function of the first kind of order zero.

\begin{figure}[htbp!]
\centering
\subfigure[Error versus CPU time.]{
\includegraphics[width=0.49\linewidth]{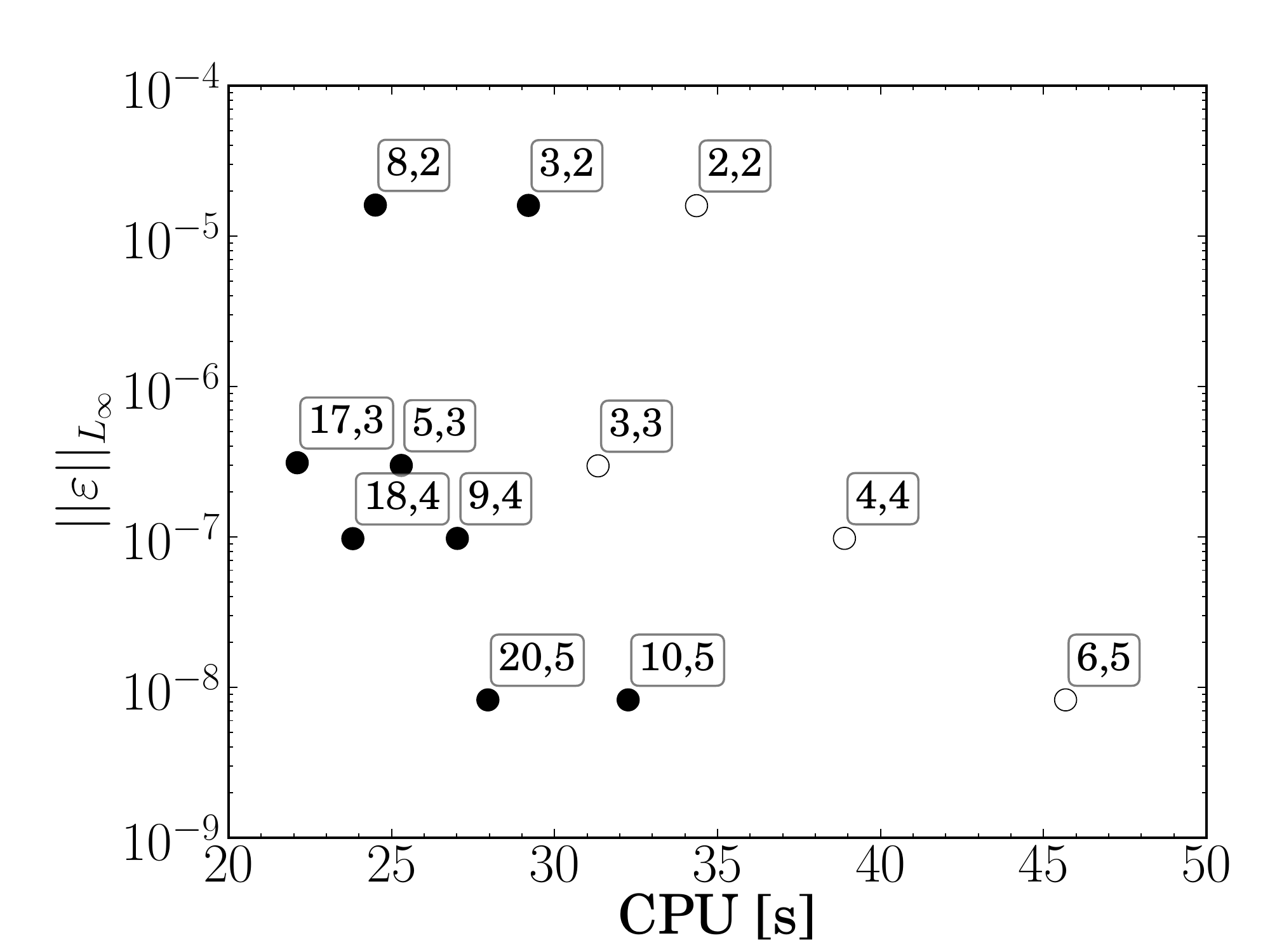}
\label{fig:error-cpu-acoustic}}
\hspace{-0.72cm}
\subfigure[CPU time for each simulation.]{
\includegraphics[width=0.49\linewidth]{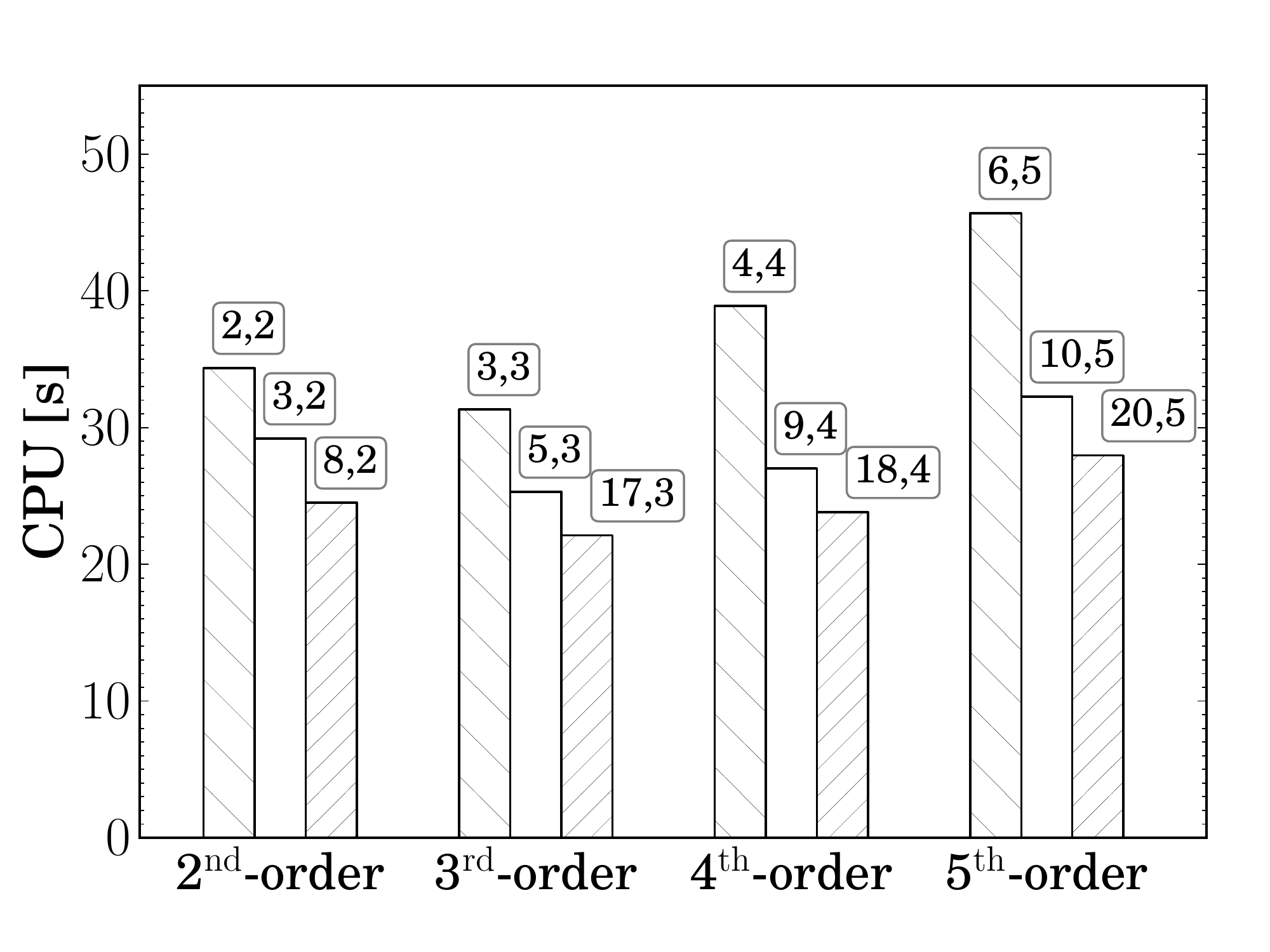}
\label{fig:cpu-acoustic}}
\caption{Error and CPU time for the acoustic pulse problem.
The labels $s,p$ for each point indicate the number of stages $s$ and the order $p$
of the corresponding scheme.  Open circles are used for the reference methods;
closed circles are used for the optimized methods.\label{fig:err-cpu-sub-acoustic}}
\end{figure}

Figure \ref{fig:err-cpu-sub-acoustic} shows the maximum norm of the error
and the CPU time for each scheme, using the predicted maximum stable CFL
number.  All schemes are again stable at their respective theoretical CFL
values.  The results are similar
to those shown in Figure \ref{fig:err-cpu-sub-annulus} for advection,
although it appears that the spatial errors are even more
dominant for this problem as the overall error is nearly the
same for all time-stepping schemes of a given order. 

\begin{figure}[htbp!]
\centering
\subfigure[Optimal 2$^\mathrm{nd}$-order methods vs. ERK(2,2).]{
\includegraphics[scale=0.375]{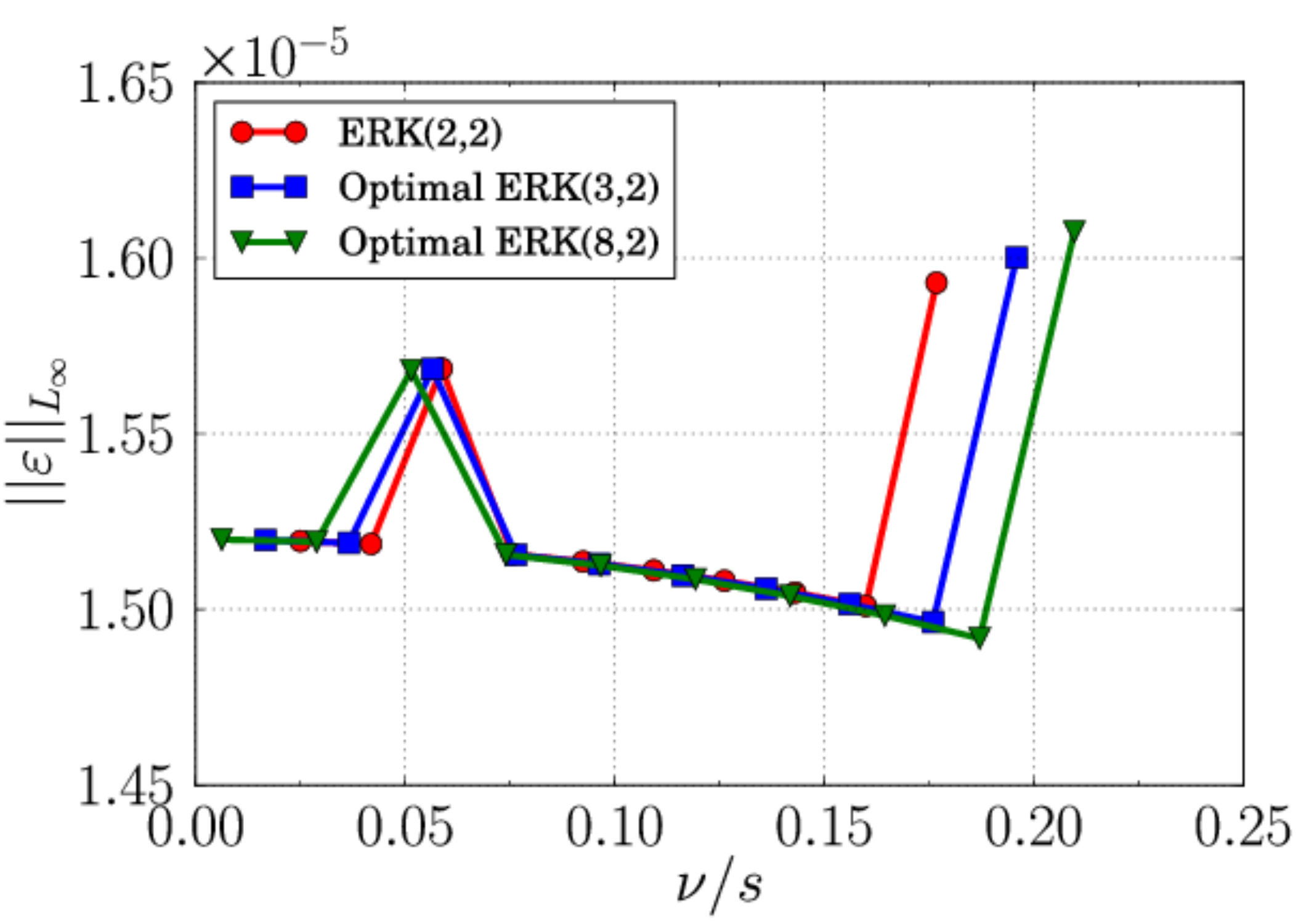}
\label{fig:acoustic_pulse-optim-erk2}}
\hspace{-0.57cm}
\subfigure[Optimal 3$^\mathrm{rd}$-order methods vs. ERK(3,3).]{
\includegraphics[scale=0.375]{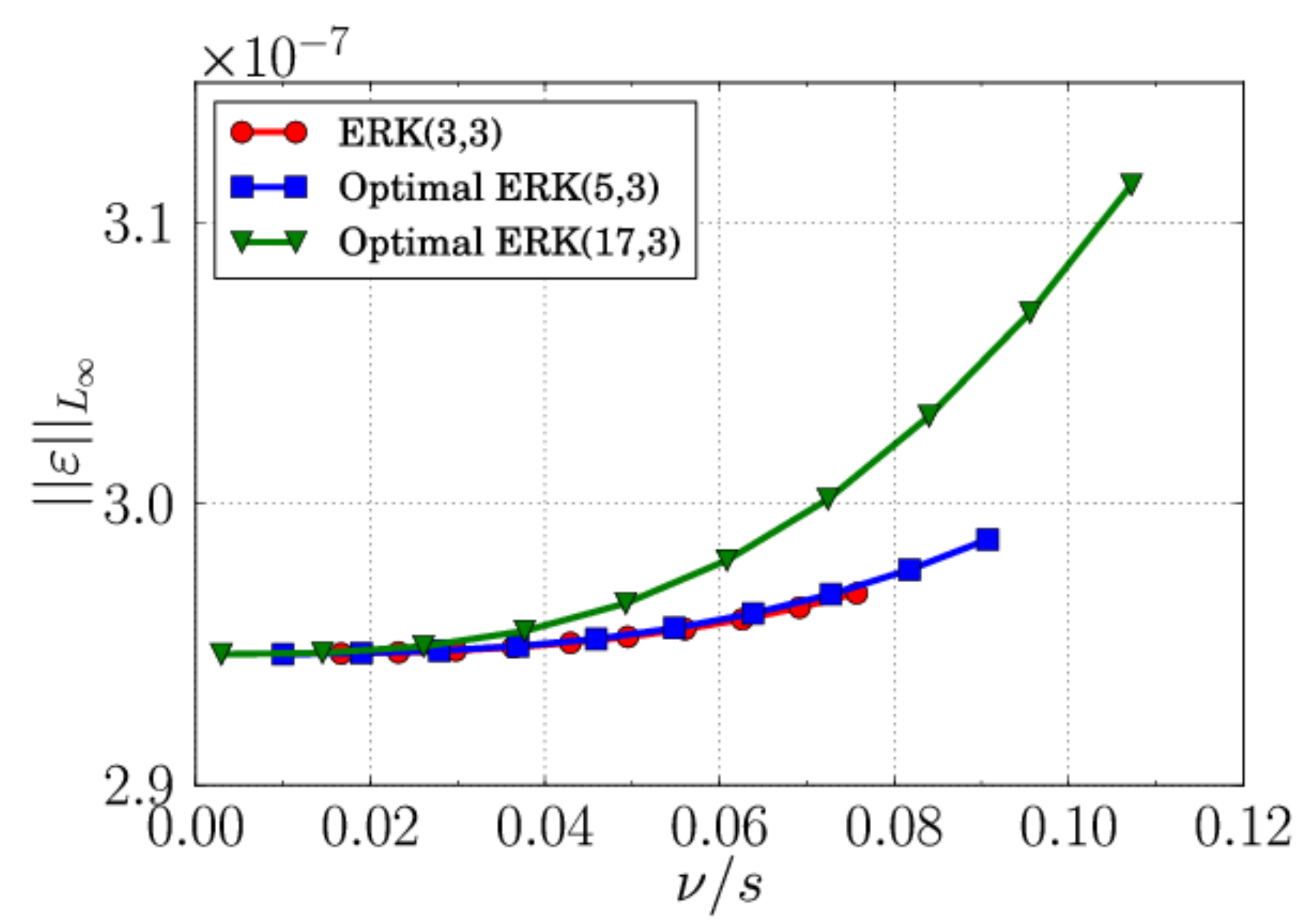}
\label{fig:acoustic_pulse-optim-erk3}}

\subfigure[Optimal 4$^\mathrm{th}$-order methods vs. ERK(4,4).]{
\includegraphics[scale=0.375]{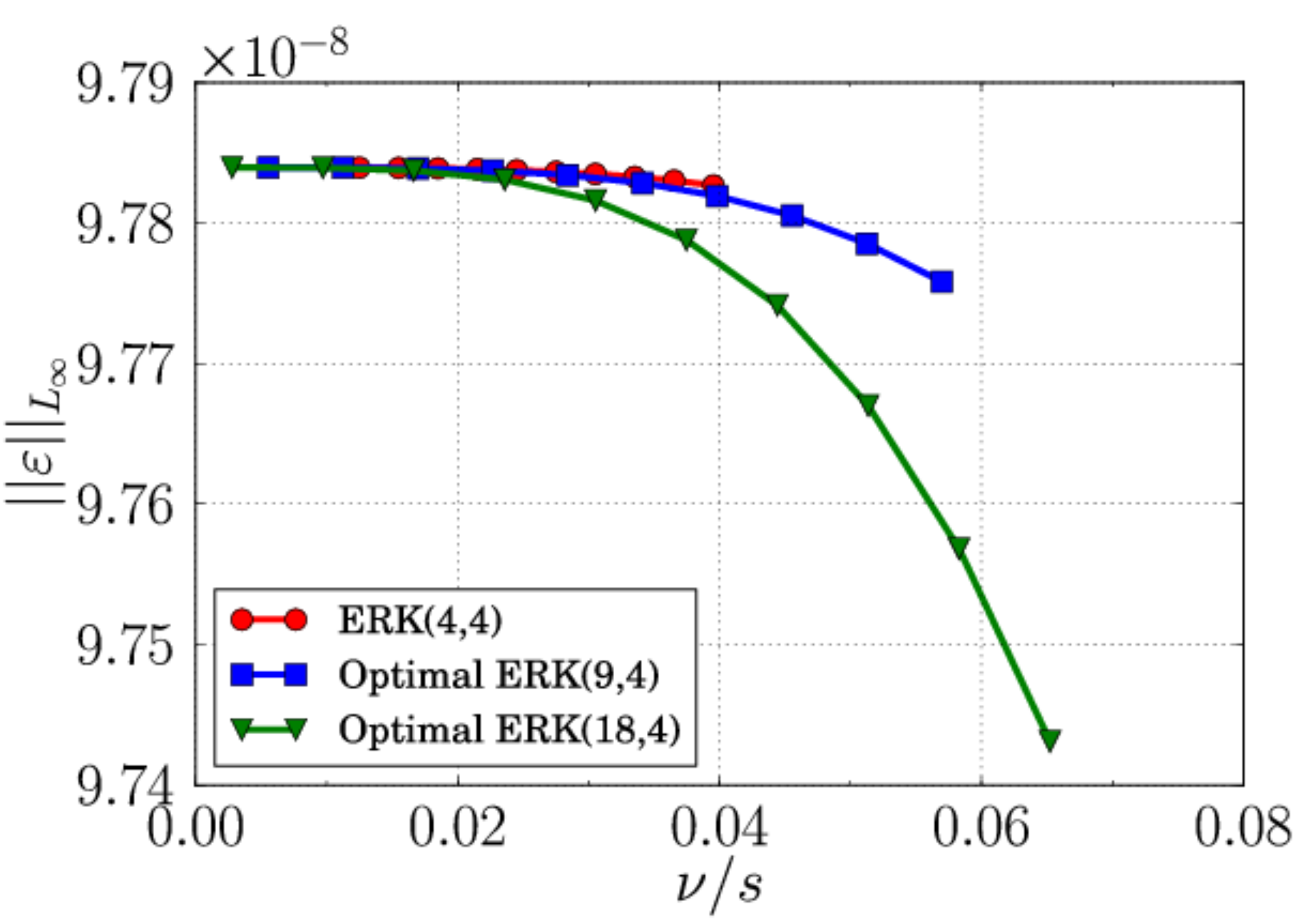}
\label{fig:acoustic_pulse-optim-erk4}}
\hspace{-0.57cm}
\subfigure[Optimal 5$^\mathrm{th}$-order methods vs. ERK(6,5).]{
\includegraphics[scale=0.375]{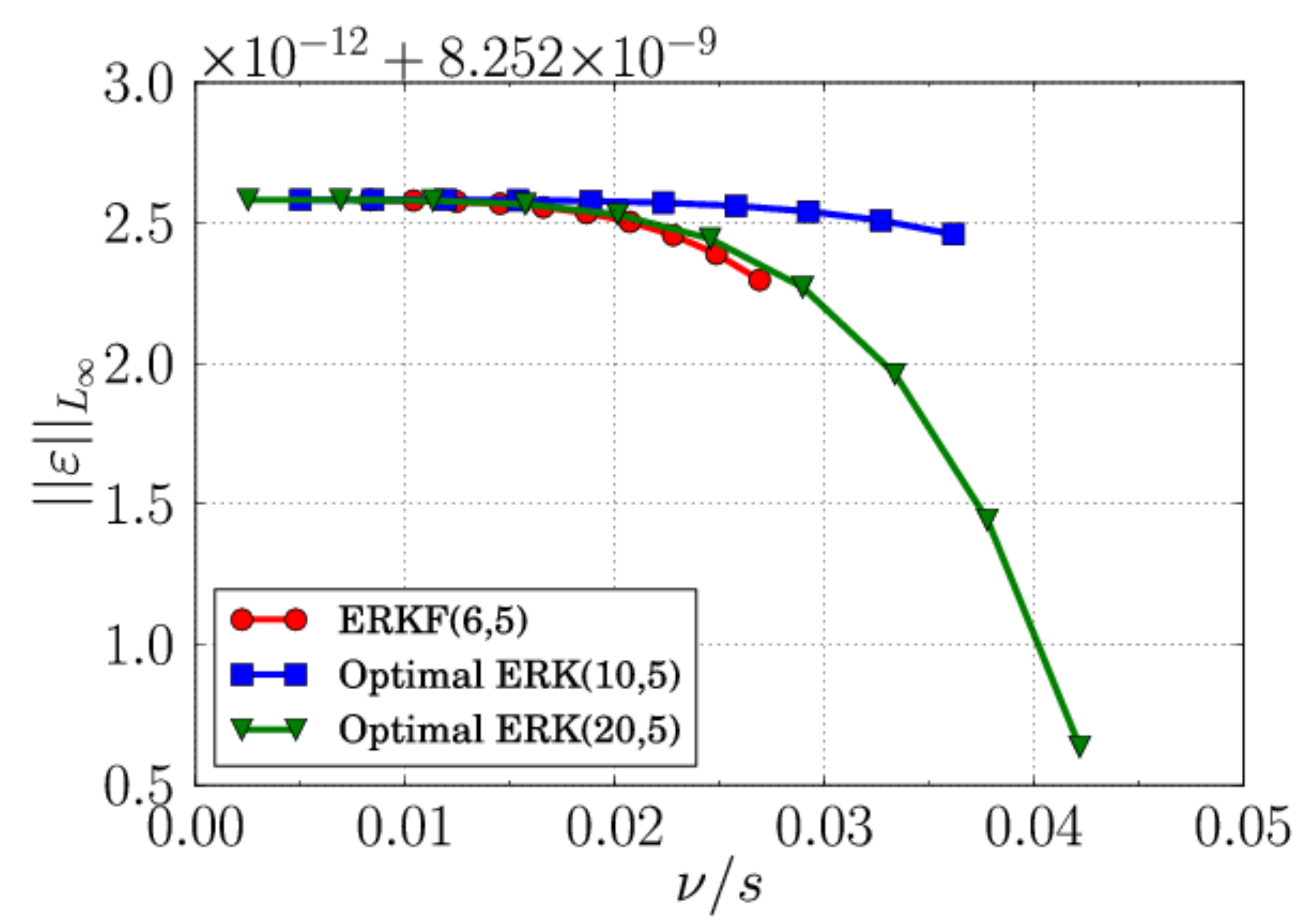}
\label{fig:acoustic_pulse-optim-erk5}}

\caption{Influence of the CFL number on the maximum error norm 
of a 2D acoustic pulse propagating in a circular domain.\label{fig:acoustic_pulse}}
\end{figure}

Figure \ref{fig:acoustic_pulse} shows the maximum norm of the error 
versus effective CFL number for a range of CFL numbers.
We observe that for a fixed order of accuracy the 
error is almost independent of the CFL number. 
More precisely we find 
that as long as the time step is smaller than the theoretical maximum stable
value, the error is dominated by the spatial discretization error. 
Note that a unit increment of the order of accuracy 
of the full discretization leads to a reduction of the error of one order of magnitude.
The optimized RK schemes speed up the simulations considerably.

\subsection{Vortex shedding past a wedge}
This test case focuses on the von Karman vortex street past a triangular wedge
\cite{Wang2007d} computed with the compressible Euler equations \eqref{eq:Euler}.
Indeed, in the inviscid framework, vortex shedding phenomena can be 
described when the considered bodies present sharp corners which ensure the 
separation of the flow.
This numerical test represents a more realistic application and it is used to 
study the performance of the new temporal schemes for a nonlinear system of 
PDEs and highly unstructured mesh. The compressible Euler equations are generally used to model the flow of an 
inviscid fluid, or the flow of a viscous fluid in regions 
where the effects of viscosity and heat conduction are negligible.
Typical applications include preliminary aircraft design and 
rotor-flow computations.

\begin{figure}
\centering
\includegraphics[scale=1.0]{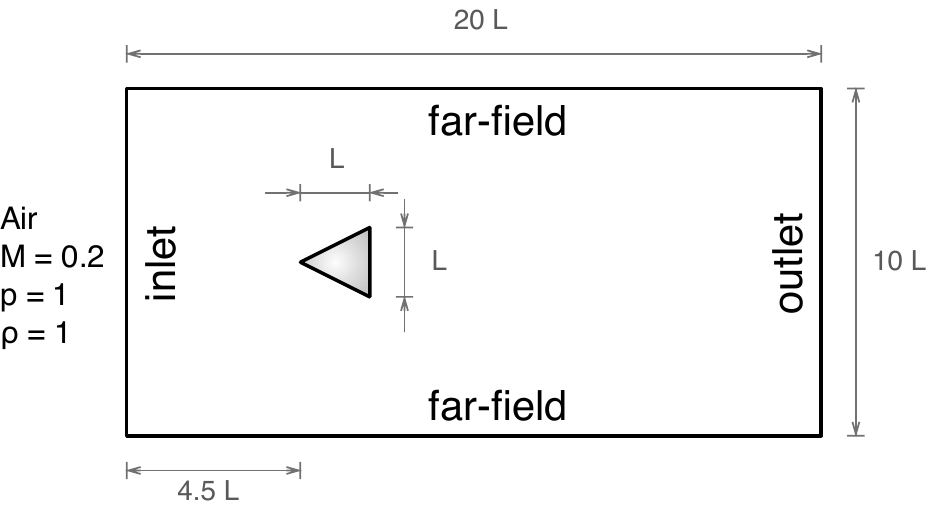}
\caption{Configuration of the wedge problem.\label{fig:configuration-wedge}}
\end{figure}

In Figure \ref{fig:configuration-wedge} the 
configuration of the test case is illustrated, where the incoming flow is from 
left to right. The wedge is placed on the centerline $y=0$ of the computational
domain and it is characterized by a length L. At the left boundary (the inflow) 
the flow is prescribed to be
uniform with zero angle of attack and free-stream Mach number of $0.2$. Both 
inlet density and inlet pressure are set to one. A pressure outlet boundary 
condition is imposed on the right boundary of the domain which is placed 
about 15 L away from the wedge. Far-field boundary conditions (i.e. uniform
Dirichlet boundary conditions for the conserved variables) are imposed both on the
top and bottom boundaries.

An unstructured grid with $11686$ quadrilateral cells with a maximum aspect ratio of 
$1.78$ and a maximum skewness value of $0.54$ is used for the 2$^\mathrm{nd}$-order calculations. This leads 
to $46744$ DOFs. The number of DOFs are again kept
about constant for higher order accurate simulation by coarsening the grid.
For this test the exact solution is not available. Therefore a reference 
solution is numerically computed by solving the problem on the mesh with $11686$ 
quadrilateral cells with the 5$^\mathrm{th}$-order SD method ($292150$ DOFs) and the ERKF(6,5) 
scheme. A CFL number $\nu = 0.1$ is used for the reference computation.

\begin{figure}[ht!]
\centering
\includegraphics[scale=0.4]{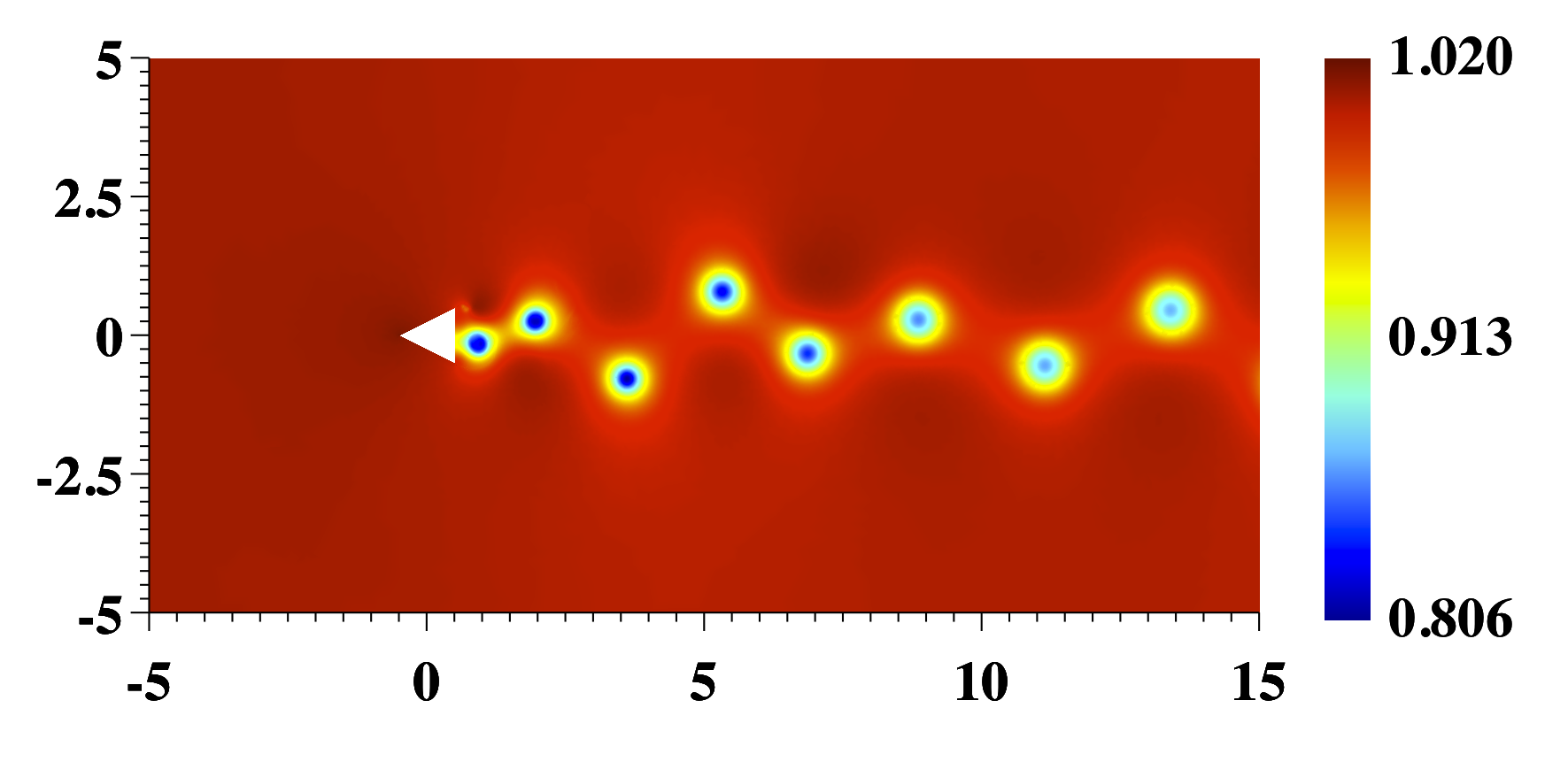}
\caption{Density contour of the flow past a wedge at $t = 200$; solution computed 
with the 4$^\mathrm{th}$-order SD method and the optimal ERK(18,4) scheme.}
\label{fig:wedge_contour_4th_order}
\end{figure}

\begin{figure}[htbp!]
\centering
\subfigure[Error versus CPU time.]{
\includegraphics[width=0.49\linewidth]{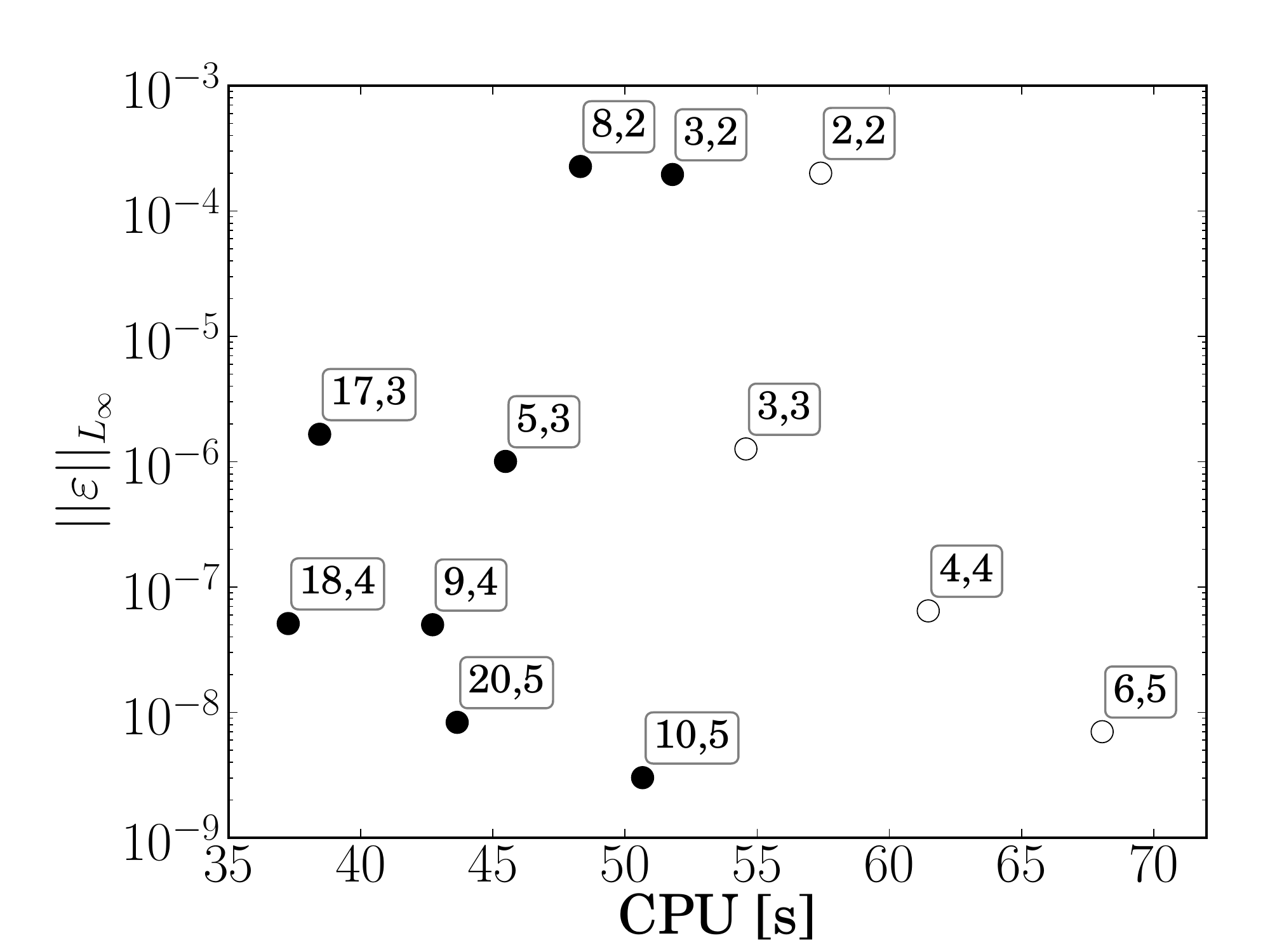}
\label{fig:error-cpu-wedge}}
\hspace{-0.72cm}
\subfigure[CPU time for each simulation.]{
\includegraphics[width=0.49\linewidth]{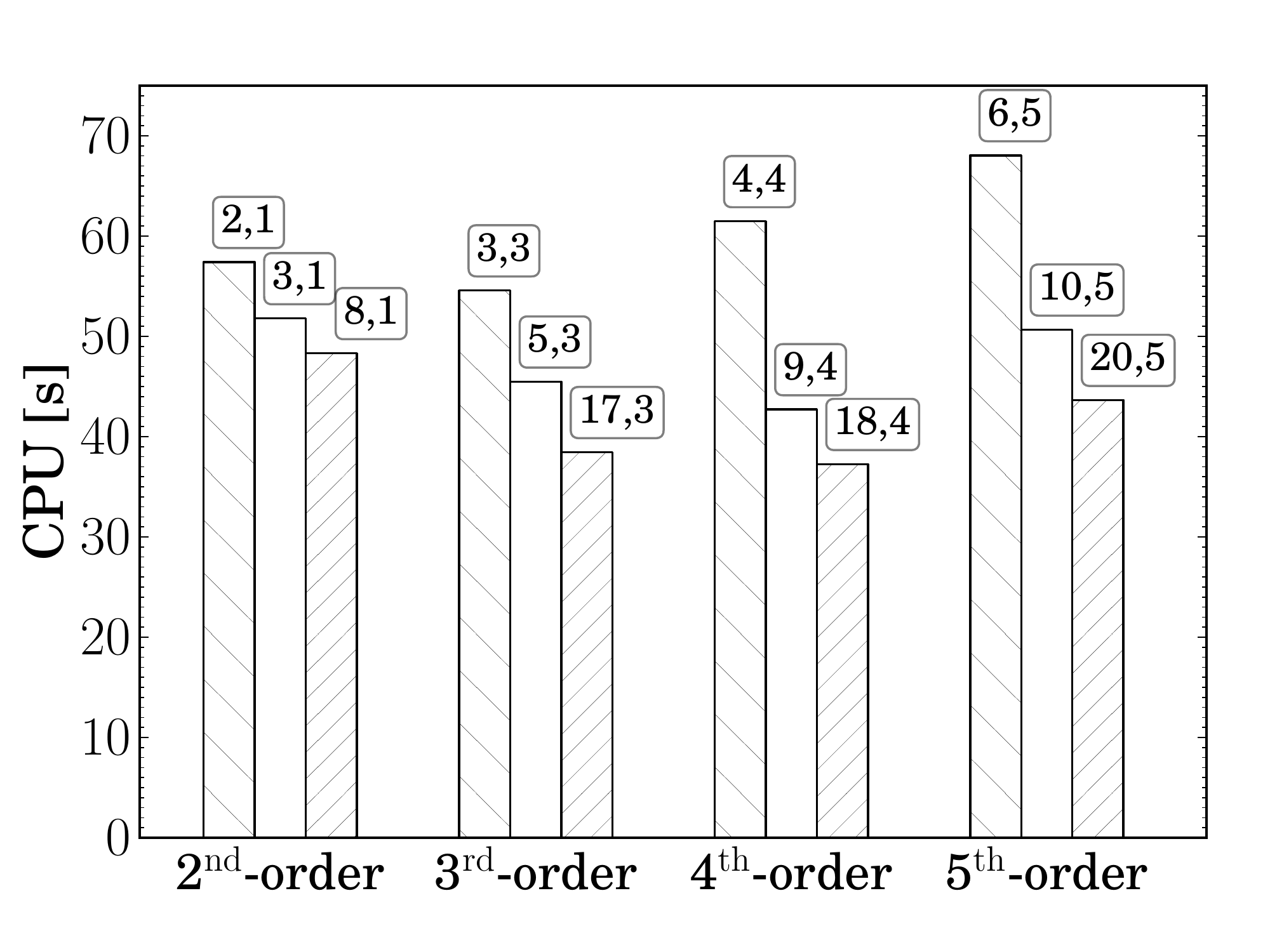}
\label{fig:cpu-wedge}}
\caption{Error and CPU time for the wedge problem.
The labels $s,p$ for each point indicate the number of stages $s$ and the order $p$
of the corresponding scheme.  Open circles are used for the reference methods;
closed circles are used for the optimized methods.\label{fig:err-cpu-sub-wedge}}
\end{figure}

In order to avoid discontinuities near the 
surface of the wedge during the transitional phase that is produced by 
the uniform free-stream initial conditions, an intermediate solution, in which the formed
vortices have not yet separated, is computed with 1$^\mathrm{st}$-order SD and 
ERK(2,2). That solution is used as the initial condition for all higher-order
computations (including the reference one) which are carried out from $t^0 = 0$
to $t^e=200$ (see Figure \ref{fig:wedge_contour_4th_order}) to generate new
unsteady laminar initial solutions for each order of accuracy. 
Afterwards, starting with these intermediate solutions,
several computations are performed using the CFL number $\dtm$ for each scheme and 
measuring the error after 0.1 seconds. Figure \ref{fig:err-cpu-sub-wedge} shows
the maximum norm of the error
and the CPU time for each scheme. 
Remarkably, we observe that the new schemes, designed using linear advection 
on a uniform grid, perform very well for the 
compressible Euler equations on an unstructured grid.
Indeed, they speed up the simulations considerably, while 
retaining a small error norm. Moreover, we highlight that the use of a
quasi-uniform grid and the CFL number $\dtm$ results in stable 
full discretizations. Therefore, also for this nonlinear test the theoretical 
stability efficiencies obtained in the first optimization step by using the 2D 
linear advection equation model are recovered. 

\section{Conclusions and future work}\label{sec:conclusions}
In this work we have developed new robust optimized explicit 
Runge-Kutta schemes for the spectral difference method
to efficiently and accurately solve wave propagation problems on unstructured 
uniform or quasi-uniform non-simplex cell grids. We have shown that by using
low-storage schemes with optimized stability function and reasonable leading truncation
error constant, one can significantly improve the performance of the resulting method of lines 
discretization. By integrating 
high-order accurate spectral difference semi-discretizations
(i.e. 3$^\mathrm{rd}$-, 4$^\mathrm{th}$- and 5$^\mathrm{th}$-order) with optimized
Runge--Kutta methods, we have found stability efficiency improvements
of $42\%$ to $65\%$, for typical systems of hyperbolic conservation laws 
used in fluid dynamics. These performance gains correspond to a reduction
in computational cost of $29\%$ to $40\%$ for a fixed simulation time. 
These improvements, which agree remarkably well with theoretical predictions
based on analysis of the 2D advection equation, are obtained without significant
sacrifices in accuracy. Indeed, when both spatial and temporal discretizations have
the same order of accuracy the spatial error typically dominates. Therefore,
for a fixed order of accuracy the error is almost independent of the CFL number
and large time steps can be used. 

Our results also highlight the advantage of high order schemes, which is even
more pronounced when optimized time integrators are employed.  
The schemes designed in this work are intended for the solution of purely
hyperbolic systems of conservation discretized on unstructured grids. The goal of ongoing research involves the optimization 
of such a family of schemes for convection-dominated problems with diffusion, where 
anisotropic grids are needed to discretize the domain with an economical 
distribution of cells.

Finally, we want to highlight that by keeping the number of
degrees of freedom constant, while increasing the order of accuracy, the new schemes allow
to get much more accurate solutions in about the same time, with about the same memory
requirements. 
For example, using our optimized
schemes, the proposed 5$^\mathrm{th}$-order discretization is actually {\em faster} than
the 2$^\mathrm{nd}$-order discretization, when the number of degrees of freedom is
held constant.  Of course, the 5$^\mathrm{th}$-order discretization is also {\em much}
more accurate.
We expect that similar improvements could be obtained by using
our approach to design optimized ERK schemes for other high-order accurate 
semi-discretizations.

\appendix



\section{Low storage ERK coefficients}\label{app:rk_coeffs}

\begin{table}[ht]
\centering
\resizebox{\linewidth}{!}{\begin{tabular}{r|r|r|r|r|r}
\multicolumn{1}{c|}{c}  &  \multicolumn{1}{c|}{$\beta$}  &  \multicolumn{1}{c|}{$\gamma_1$}  &  \multicolumn{1}{c|}{$\gamma_2$}  &  \multicolumn{1}{c|}{$\gamma_3$} &  \multicolumn{1}{c}{$\delta$} \\
\hline
 0.0000000000000000\e{+00}  &  7.2366074728360086\e{-01}  &  0.0000000000000000\e{+00}  &  1.0000000000000000\e{+00}  &  0.0000000000000000\e{+00} &  1.0000000000000000\e{+00} \\
 7.2366074728360086\e{-01}  &  3.4217876502651023\e{-01}  &  -1.2664395576322218\e{-01}  &  6.5427782599406470\e{-01}  &  0.0000000000000000\e{+00} &  7.2196567116037724\e{-01} \\
 5.9236433182015646\e{-01}  &  3.6640216242653251\e{-01}  &  1.1426980685848858\e{+00}  &  -8.2869287683723744\e{-02}  &  0.0000000000000000\e{+00} &  0.0000000000000000\e{+00} \\
\end{tabular}}
\caption{3S* low storage coefficients of ERK(3,2) method, optimized for 2$^\mathrm{nd}$-order SD scheme.}
\label{tbl:rk-coeffs-3Sstar-erk2-3}
\end{table} 

\begin{table}[ht]
\centering
\resizebox{\linewidth}{!}{\begin{tabular}{r|r|r|r|r|r}
\multicolumn{1}{c|}{c}  &  \multicolumn{1}{c|}{$\beta$}  &  \multicolumn{1}{c|}{$\gamma_1$}  &  \multicolumn{1}{c|}{$\gamma_2$}  &  \multicolumn{1}{c|}{$\gamma_3$} &  \multicolumn{1}{c}{$\delta$} \\
\hline
 0.0000000000000000\e{+00}  &  9.9292229393265474\e{-01}  &  0.0000000000000000\e{+00}  &  1.0000000000000000\e{+00}  &  0.0000000000000000\e{+00} &  1.0000000000000000\e{+00} \\
 9.9292229393265474\e{-01}  &  5.2108385130005974\e{-01}  &  4.2397552118208004\e{-01}  &  4.4390665802303775\e{-01}  &  0.0000000000000000\e{+00} &  2.9762522910396538\e{-01} \\
 1.0732413280565014\e{+00}  &  3.8505327083543915\e{-03}  &  -2.3528852074619033\e{-01}  &  7.5333732286056154\e{-01}  &  0.0000000000000000\e{+00} &  3.4212961014330662\e{-01} \\
 2.5057060509809409\e{-01}  &  7.9714199213087467\e{-01}  &  7.9598685017877846\e{-01}  &  6.5885460813015481\e{-02}  &  5.8415358412023582\e{-02} &  5.7010739154759105\e{-01} \\
 1.0496674928979783\e{+00}  &  -8.1822460276649120\e{-02}  &  -1.3205224623823271\e{+00}  &  6.3976199384289623\e{-01}  &  6.4219008773865116\e{-01} &  4.1350769551529132\e{-01} \\
 -6.7488037049720317\e{-01}  &  8.4604310411858186\e{-01}  &  2.1452956294251941\e{+00}  &  -7.3823030755143193\e{-01}  &  6.8770305706885126\e{-01} &  -1.4040672669058066\e{-01} \\
 -1.5868411612120166\e{+00}  &  -1.0191166090841246\e{-01}  &  -9.5532770501880648\e{-01}  &  7.0177211879534529\e{-01}  &  6.3729822311671305\e{-02} &  2.1249567092409008\e{-01} \\
 2.1138242369563969\e{+00}  &  6.3190236038107500\e{-02}  &  2.5361391125131094\e{-01}  &  4.0185379950224559\e{-01}  &  -3.3679429978131387\e{-01} &  0.0000000000000000\e{+00} \\ 
\end{tabular}}
\caption{3S* low storage coefficients of ERK(8,2) method, optimized for 2$^\mathrm{nd}$-order SD scheme.}
\label{tbl:rk-coeffs-3Sstar-erk2-8}
\end{table} 

\begin{table}[ht]
\centering
\resizebox{\linewidth}{!}{
\begin{tabular}{r|r|r|r|r|r}
\multicolumn{1}{c|}{c}  &  \multicolumn{1}{c|}{$\beta$}  &  \multicolumn{1}{c|}{$\gamma_1$}  &  \multicolumn{1}{c|}{$\gamma_2$}  &  \multicolumn{1}{c|}{$\gamma_3$} &  \multicolumn{1}{c}{$\delta$} \\
\hline
 0.0000000000000000\e{+00}  &  2.3002859824852059\e{-01}  &  0.0000000000000000\e{+00}  &  1.0000000000000000\e{+00}  &  0.0000000000000000\e{+00} &  1.0000000000000000\e{+00} \\
 2.3002859824852059\e{-01}  &  3.0214498165167158\e{-01}  &  2.5876919610938998\e{-01}  &  5.5284013909611196\e{-01}  &  0.0000000000000000\e{+00} &  3.4076878915216791\e{-01} \\
 4.0500453764839639\e{-01}  &  8.0256010238856679\e{-01}  &  -1.3243708384977859\e{-01}  &  6.7318513326032769\e{-01}  &  0.0000000000000000\e{+00} &  3.4143871647890728\e{-01} \\
 8.9478204142351003\e{-01}  &  4.3621618871511753\e{-01}  &  5.0556648948362981\e{-02}  &  2.8031054965521607\e{-01}  &  2.7525797946334213\e{-01} &  7.2292984084963252\e{-01} \\
 7.2351146275625733\e{-01}  &  1.1292705979513513\e{-01}  &  5.6705507883024708\e{-01}  &  5.5215115815918758\e{-01}  &  -8.9505445022148511\e{-01} &  0.0000000000000000\e{+00} \\
\end{tabular}}
\caption{3S* low storage coefficients of ERK(5,3) method, optimized for 3$^\mathrm{rd}$-order SD scheme.}
\label{tbl:rk-coeffs-3Sstar-erk3-5}
\end{table} 

\begin{table}[ht]
\centering
\resizebox{\linewidth}{!}{\begin{tabular}{r|r|r|r|r|r}
\multicolumn{1}{c|}{c}  &  \multicolumn{1}{c|}{$\beta$}  &  \multicolumn{1}{c|}{$\gamma_1$}  &  \multicolumn{1}{c|}{$\gamma_2$}  &  \multicolumn{1}{c|}{$\gamma_3$} &  \multicolumn{1}{c}{$\delta$} \\
\hline
 0.0000000000000000\e{+00}  &  4.9565403010221741\e{-02}  &  0.0000000000000000\e{+00}  &  1.0000000000000000\e{+00}  &  0.0000000000000000\e{+00} &  1.0000000000000000\e{+00} \\
 4.9565403010221741\e{-02}  &  9.7408718698159397\e{-02}  &  7.9377023961829174\e{-01}  &  3.2857861940811250\e{-01}  &  0.0000000000000000\e{+00} &  -3.7235794357769936\e{-01} \\
 1.3068799001687578\e{-01}  &  -1.7620737976801870\e{-01}  &  -8.3475116244241754\e{-02}  &  1.1276843361180819\e{+00}  &  0.0000000000000000\e{+00} &  3.3315440189685536\e{-01} \\
 -1.5883063460310493\e{-01}  &  1.4852069175460250\e{-01}  &  -1.6706337980062214\e{-02}  &  1.3149447395238016\e{+00}  &  8.4034574578399479\e{-01} &  -8.2667630338402520\e{-01} \\
 3.5681144740196935\e{-01}  &  -3.3127657103714951\e{-02}  &  3.6410691500331427\e{-01}  &  5.2062891534209055\e{-01}  &  8.5047738439705145\e{-01} &  -5.4628377681035534\e{-01} \\
 7.6727123317642698\e{-02}  &  4.8294609330498492\e{-02}  &  6.9178255181542780\e{-01}  &  8.8127462325164985\e{-01}  &  1.4082448501410852\e{-01} &  6.0210777634642887\e{-01} \\
 1.0812579255374613\e{-01}  &  4.9622612199980112\e{-02}  &  1.4887115004739182\e{+00}  &  4.2020606445856712\e{-01}  &  -3.2678802469519369\e{-01} &  -5.7528717894031067\e{-01} \\
 1.8767228084815801\e{-01}  &  8.7340766269850378\e{-01}  &  4.5336125560871188\e{-01}  &  7.6532635739246124\e{-02}  &  5.3716357620635535\e{-01} &  5.0914861529202782\e{-01} \\
 9.6162976936182631\e{-01}  &  -2.8692804399085370\e{-01}  &  -1.2705776046458739\e{-01}  &  4.4386734924685722\e{-01}  &  9.0228922115199051\e{-01} &  3.8258114767897194\e{-01} \\
 -2.2760719867560897\e{-01}  &  1.2679897532256112\e{+00}  &  8.3749845457747696\e{-01}  &  6.6503093955199682\e{-02}  &  1.5960226946983552\e{-01} &  -4.6279063221185290\e{-01} \\
 1.1115681606027146\e{+00}  &  -1.0217436118953449\e{-02}  &  1.5709218393361746\e{-01}  &  1.5850209163184039\e{+00}  &  1.1038153140686748\e{+00} &  -2.0820434288562648\e{-01} \\
 6.1266845427676520\e{-01}  &  8.4665570032598350\e{-02}  &  -5.7768207086288348\e{-01}  &  1.1521721573462576\e{+00}  &  1.0843516423068365\e{-01} &  1.4398056081552713\e{+00} \\
 1.0729473245077408\e{+00}  &  2.8253854742588246\e{-02}  &  -5.7340394122375393\e{-01}  &  1.1172750819374575\e{+00}  &  4.6212710442787724\e{-01} &  -2.8056600927348752\e{-01} \\
 3.7824186468104548\e{-01}  &  -9.2936733010804407\e{-02}  &  -1.2050734846514470\e{+00}  &  7.7630223917584007\e{-01}  &  -3.3448312125108398\e{-01} &  2.2767189929551406\e{+00} \\
 7.9041891347646720\e{-01}  &  -8.4798124766803512\e{-02}  &  -2.8100719513641002\e{+00}  &  1.0046657060652295\e{+00}  &  1.1153826567096696\e{+00} &  -5.8917530100546356\e{-01} \\
 -1.0406955693161675\e{+00}  &  -1.6923145636158564\e{-02}  &  1.6142798657609492\e{-01}  &  -1.9795868964959054\e{-01}  &  1.5503248734613539\e{+00} &  9.1328651048418164\e{-01} \\
 -2.4607146824557105\e{-01}  &  -4.7305106233879957\e{-02}  &  -2.5801264756641613\e{+00}  &  1.3350583594705518\e{+00}  &  -1.2200245424704212\e{+00} &  0.0000000000000000\e{+00} \\
\end{tabular}}
\caption{3S* low storage coefficients of ERK(17,3) method, optimized for 3$^\mathrm{rd}$-order SD scheme.}
\label{tbl:rk-coeffs-3Sstar-erk3-17}
\end{table} 

\begin{table}[ht]
\centering
\resizebox{\linewidth}{!}{\begin{tabular}{r|r|r|r|r|r}
\multicolumn{1}{c|}{c}  &  \multicolumn{1}{c|}{$\beta$}  &  \multicolumn{1}{c|}{$\gamma_1$}  &  \multicolumn{1}{c|}{$\gamma_2$}  &  \multicolumn{1}{c|}{$\gamma_3$} &  \multicolumn{1}{c}{$\delta$} \\
\hline
 0.0000000000000000\e{+00}  &  2.8363432481011769\e{-01}  &  0.0000000000000000\e{+00}  &  1.0000000000000000\e{+00}  &  0.0000000000000000\e{+00} &  1.0000000000000000\e{+00} \\
 2.8363432481011769\e{-01}  &  9.7364980747486463\e{-01}  &  -4.6556413837561301\e{+00}  &  2.4992627683300688\e{+00}  &  0.0000000000000000\e{+00} &  1.2629238731608268\e{+00} \\
 5.4840742446661772\e{-01}  &  3.3823592364196498\e{-01}  &  -7.7202649689034453\e{-01}  &  5.8668202764174726\e{-01}  &  0.0000000000000000\e{+00} &  7.5749675232391733\e{-01} \\
 3.6872298094969475\e{-01}  &  -3.5849518935750763\e{-01}  &  -4.0244202720632174\e{+00}  &  1.2051419816240785\e{+00}  &  7.6209857891449362\e{-01} &  5.1635907196195419\e{-01} \\
 -6.8061183026103156\e{-01}  &  -4.1139587569859462\e{-03}  &  -2.1296873883702272\e{-02}  &  3.4747937498564541\e{-01}  &  -1.9811817832965520\e{-01} &  -2.7463346616574083\e{-02} \\
 3.5185265855105619\e{-01}  &  1.4279689871485013\e{+00}  &  -2.4350219407769953\e{+00}  &  1.3213458736302766\e{+00}  &  -6.2289587091629484\e{-01} &  -4.3826743572318672\e{-01} \\
 1.6659419385562171\e{+00}  &  1.8084680519536503\e{-02}  &  1.9856336960249132\e{-02}  &  3.1196363453264964\e{-01}  &  -3.7522475499063573\e{-01} &  1.2735870231839268\e{+00} \\
 9.7152778807463247\e{-01}  &  1.6057708856060501\e{-01}  &  -2.8107894116913812\e{-01}  &  4.3514189245414447\e{-01}  &  -3.3554373281046146\e{-01} &  -6.2947382217730230\e{-01} \\
 9.0515694340066954\e{-01}  &  2.9522267863254809\e{-01}  &  1.6894354373677900\e{-01}  &  2.3596980658341213\e{-01}  &  -4.5609629702116454\e{-02} &  0.0000000000000000\e{+00} \\
\end{tabular}}
\caption{3S* low storage coefficients of ERK(9,4) method, optimized for 4$^\mathrm{th}$-order SD scheme.}
\label{tbl:rk-coeffs-3Sstar-erk4-9}
\end{table} 

\begin{table}[ht]
\centering
\resizebox{\linewidth}{!}{\begin{tabular}{r|r|r|r|r|r}
\multicolumn{1}{c|}{c}  &  \multicolumn{1}{c|}{$\beta$}  &  \multicolumn{1}{c|}{$\gamma_1$}  &  \multicolumn{1}{c|}{$\gamma_2$}  &  \multicolumn{1}{c|}{$\gamma_3$} &  \multicolumn{1}{c}{$\delta$} \\
\hline
 0.0000000000000000\e{+00}  &  1.2384169480626298\e{-01}  &  0.0000000000000000\e{+00}  &  1.0000000000000000\e{+00}  &  0.0000000000000000\e{+00} &  1.0000000000000000\e{+00} \\
 1.2384169480626298\e{-01}  &  1.0176262534280349\e{+00}  &  1.1750819811951678\e{+00}  &  -1.2891068509748144\e{-01}  &  0.0000000000000000\e{+00} &  3.5816500441970289\e{-01} \\
 1.1574324659554065\e{+00}  &  -6.9732026387527429\e{-02}  &  3.0909017892654811\e{-01}  &  3.5609406666728954\e{-01}  &  0.0000000000000000\e{+00} &  5.8208024465093577\e{-01} \\
 5.4372099141546926\e{-01}  &  3.4239356067806476\e{-01}  &  1.4409117788115862\e{+00}  &  -4.0648075226104241\e{-01}  &  2.5583378537249163\e{-01} &  -2.2615285894283538\e{-01} \\
 8.8394666834280744\e{-01}  &  1.8177707207807942\e{-02}  &  -4.3563049445694069\e{-01}  &  6.0714786995207426\e{-01}  &  5.2676794366988289\e{-01} &  -2.1715466578266213\e{-01} \\
 -1.2212042176605774\e{-01}  &  -6.1188746289480445\e{-03}  &  2.0341503014683893\e{-01}  &  1.0253501186236846\e{+00}  &  -2.5648375621792202\e{-01} &  -4.6990441450888265\e{-01} \\
 4.4125685133082082\e{-01}  &  7.8242308902580354\e{-02}  &  4.9828356971917692\e{-01}  &  2.4411240760769423\e{-01}  &  3.1932438003236391\e{-01} &  -2.7986911594744995\e{-01} \\
 3.8039092095473748\e{-01}  &  -3.7642864750532951\e{-01}  &  3.5307737157745489\e{+00}  &  -1.2813606970134104\e{+00}  &  -3.1106815010852862\e{-01} &  9.8513926355272197\e{-01} \\
 5.4591107347528367\e{-02}  &  -4.5078383666690258\e{-02}  &  -7.9318790975894626\e{-01}  &  8.1625711892373898\e{-01}  &  4.7631196164025996\e{-01} &  -1.1899324232814899\e{-01} \\
 4.8731855535356028\e{-01}  &  -7.5734228201432585\e{-01}  &  8.9120513355345166\e{-01}  &  1.0171269354643386\e{-01}  &  -9.8853727938895783\e{-02} &  4.2821073124370562\e{-01} \\
 -2.3007964303896034\e{-01}  &  -2.7149222760935121\e{-01}  &  5.7091009196320974\e{-01}  &  1.9379378662711269\e{-01}  &  1.9274726276883622\e{-01} &  -8.2196355299900403\e{-01} \\
 -1.8907656662915873\e{-01}  &  1.1833684341657344\e{-03}  &  1.6912188575015419\e{-02}  &  7.4408643544851782\e{-01}  &  3.2389860855971508\e{-02} &  5.8113997057675074\e{-02} \\
 8.1059805668623763\e{-01}  &  2.8858319979308041\e{-02}  &  1.0077912519329719\e{+00}  &  -1.2591764563430008\e{-01}  &  7.5923980038397509\e{-02} &  -6.1283024325436919\e{-01} \\
 7.7080875997868803\e{-01}  &  4.6005267586974657\e{-01}  &  -6.8532953752099512\e{-01}  &  1.1996463179654226\e{+00}  &  2.0635456088664017\e{-01} &  5.6800136190634054\e{-01} \\
 1.1712158507200179\e{+00}  &  1.8014887068775631\e{-02}  &  1.0488165551884063\e{+00}  &  4.5772068865370406\e{-02}  &  -8.9741032556032857\e{-02} &  -3.3874970570335106\e{-01} \\
 1.2755351018003545\e{+00}  &  -1.5508175395461857\e{-02}  &  8.3647761371829943\e{-01}  &  8.3622292077033844\e{-01}  &  2.6899932505676190\e{-02} &  -7.3071238125137772\e{-01} \\
 8.0422507946168564\e{-01}  &  -4.0095737929274988\e{-01}  &  1.3087909830445710\e{+00}  &  -1.4179124272450148\e{+00}  &  4.1882069379552307\e{-02} &  8.3936016960374532\e{-02} \\
 9.7508680250761848\e{-01}  &  1.4949678367038011\e{-01}  &  9.0419681700177323\e{-01}  &  1.3661459065331649\e{-01}  &  6.2016148912381761\e{-02} &  0.0000000000000000\e{+00} \\
\end{tabular}}
\caption{3S* low storage coefficients of ERK(18,4) method, optimized for 4$^\mathrm{th}$-order SD scheme.}
\label{tbl:rk-coeffs-3Sstar-erk4-18}
\end{table} 

\begin{table}[ht]
\centering
\resizebox{\linewidth}{!}{\begin{tabular}{r|r|r|r|r|r}
\multicolumn{1}{c|}{c}  &  \multicolumn{1}{c|}{$\beta$}  &  \multicolumn{1}{c|}{$\gamma_1$}  &  \multicolumn{1}{c|}{$\gamma_2$}  &  \multicolumn{1}{c|}{$\gamma_3$} &  \multicolumn{1}{c}{$\delta$} \\
\hline
 0.0000000000000000\e{+00}  &  2.5978835757039448\e{-01}  &  0.0000000000000000\e{+00}  &  1.0000000000000000\e{+00}  &  0.0000000000000000\e{+00} &  1.0000000000000000\e{+00} \\
 2.5978835757039448\e{-01}  &  1.7770088002098183\e{-02}  &  4.0436600785287713\e{-01}  &  6.8714670697294733\e{-01}  &  0.0000000000000000\e{+00} &  -1.3317784091400336\e{-01} \\
 9.9045731158085557\e{-02}  &  2.4816366373161344\e{-01}  &  -8.5034274641295027\e{-01}  &  1.0930247604585732\e{+00}  &  0.0000000000000000\e{+00} &  8.2604227852898304\e{-01} \\
 2.1555118823045644\e{-01}  &  7.9417368275785671\e{-01}  &  -6.9508941671218478\e{+00}  &  3.2259753823377983\e{+00}  &  -2.3934051593398129\e{+00} &  1.5137004305165804\e{+00} \\
 5.0079500784155040\e{-01}  &  3.8853912968701337\e{-01}  &  9.2387652252320684\e{-01}  &  1.0411537008416110\e{+00}  &  -1.9028544220991284\e{+00} &  -1.3058100631721905\e{+00} \\
 5.5922519148547800\e{-01}  &  1.4550516642704694\e{-01}  &  -2.5631780399589106\e{+00}  &  1.2928214888638039\e{+00}  &  -2.8200422105835639\e{+00} &  3.0366787893355149\e{+00} \\
 5.4499869734044426\e{-01}  &  1.5875173794655811\e{-01}  &  2.5457448699988827\e{-01}  &  7.3914627692888835\e{-01}  &  -1.8326984641282289\e{+00} &  -1.4494582670831953\e{+00} \\
 7.6152246625852738\e{-01}  &  1.6506056315937651\e{-01}  &  3.1258317336761454\e{-01}  &  1.2391292570651462\e{-01}  &  -2.1990945108072310\e{-01} &  3.8343138733685103\e{+00} \\
 8.4270620830633836\e{-01}  &  2.1180932999328042\e{-01}  &  -7.0071148003175443\e{-01}  &  1.8427534793568445\e{-01}  &  -4.0824306603783045\e{-01} &  4.1222939718018692\e{+00} \\
 9.1522098071770008\e{-01}  &  1.5593923403495016\e{-01}  &  4.8396209710057070\e{-01}  &  5.7127889427161162\e{-02}  &  -1.3776697911236280\e{-01} &  0.0000000000000000\e{+00} \\
\end{tabular}}
\caption{3S* low storage coefficients of ERK(10,5) method, optimized for 5$^\mathrm{th}$-order SD scheme.}
\label{tbl:rk-coeffs-3Sstar-erk5-10}
\end{table} 

\begin{table}[ht]
\centering
\resizebox{\linewidth}{!}{\begin{tabular}{r|r|r|r|r|r}
\multicolumn{1}{c|}{c}  &  \multicolumn{1}{c|}{$\beta$}  &  \multicolumn{1}{c|}{$\gamma_1$}  &  \multicolumn{1}{c|}{$\gamma_2$}  &  \multicolumn{1}{c|}{$\gamma_3$} &  \multicolumn{1}{c}{$\delta$} \\
\hline
 0.0000000000000000\e{+00}  &  1.7342385375780556\e{-01}  &  0.0000000000000000\e{+00}  &  1.0000000000000000\e{+00}  &  0.0000000000000000\e{+00} &  1.0000000000000000\e{+00} \\
 1.7342385375780556\e{-01}  &  2.8569004728564801\e{-01}  &  -1.1682479703229380\e{+00}  &  8.8952052154583572\e{-01}  &  0.0000000000000000\e{+00} &  1.4375468781258596\e{+00} \\
 3.0484982420032158\e{-01}  &  6.8727044379779589\e{-01}  &  -2.5112155037089772\e{+00}  &  8.8988129100385194\e{-01}  &  0.0000000000000000\e{+00} &  1.5081653637261594\e{+00} \\
 5.5271395645729193\e{-01}  &  1.2812121060977319\e{-01}  &  -5.5259960154735988\e{-01}  &  3.5701564494677057\e{-01}  &  1.9595487007932735\e{-01} &  -1.4575347066062688\e{-01} \\
 4.7079204549750037\e{-02}  &  4.9137180740403122\e{-04}  &  2.9243033509511740\e{-03}  &  2.4232462479216824\e{-01}  &  -6.9871675039100595\e{-05} &  3.1495761082838158\e{-01} \\
 1.5652540451324129\e{-01}  &  4.7033584446956857\e{-02}  &  -4.7948973385386493\e{+00}  &  1.2727083024258155\e{+00}  &  1.0592231169810050\e{-01} &  3.5505919368536931\e{-01} \\
 1.8602224049074517\e{-01}  &  4.4539998128170821\e{-01}  &  -5.3095533497183016\e{+00}  &  1.1126977210342681\e{+00}  &  1.0730426871909635\e{+00} &  2.3616389374566960\e{-01} \\
 2.8426620035751449\e{-01}  &  1.2259824887343720\e{+00}  &  -2.3624194456630736\e{+00}  &  5.1360709645409097\e{-01}  &  8.9257826744389124\e{-01} &  1.0267488547302055\e{-01} \\
 9.5094727548792268\e{-01}  &  2.0616463985024421\e{-02}  &  2.0068995756589547\e{-01}  &  1.1181089682044856\e{-01}  &  -1.4078912484894415\e{-01} &  3.5991243524519438\e{+00} \\
 6.8046501070096010\e{-01}  &  1.5941162575324802\e{-01}  &  -1.4985808661597710\e{+00}  &  2.7881272382085232\e{-01}  &  -2.6869890558434262\e{-01} &  1.5172890003890782\e{+00} \\
 5.9705366562360063\e{-01}  &  1.2953803678226099\e{+00}  &  4.8941228502377687\e{-01}  &  4.9032886260666715\e{-02}  &  -6.5175753568318007\e{-02} &  1.8171662741779953\e{+00} \\
 1.8970821645077285\e{+00}  &  1.7287352967302603\e{-03}  &  -1.0387512755259576\e{-01}  &  4.1871051065897870\e{-02}  &  4.9177812903108553\e{-01} &  2.8762263521436831\e{+00} \\
 2.9742664004529606\e{-01}  &  1.1660483420536467\e{-01}  &  -1.3287664273288191\e{-01}  &  4.4602463796686219\e{-02}  &  4.6017684776493678\e{-01} &  4.6350154228218754\e{-01} \\
 6.0813463700134940\e{-01}  &  7.7997036621815521\e{-02}  &  7.5858678822837511\e{-01}  &  1.4897271251154750\e{-02}  &  -6.4689512947008251\e{-03} &  1.5573122110727220\e{+00} \\
 7.3080004188477765\e{-01}  &  3.2563250234418012\e{-01}  &  -4.3321586294096939\e{+00}  &  2.6244269699436817\e{-01}  &  4.4034728024115377\e{-01} &  2.0001066778080254\e{+00} \\
 9.1656999044951792\e{-01}  &  1.0611520488333197\e{+00}  &  4.8199700138402146\e{-01}  &  -4.7486056986590294\e{-03}  &  6.1086885767527943\e{-01} &  9.1690694855534305\e{-01} \\
 1.4309687554614530\e{+00}  &  6.5891625628040993\e{-04}  &  -7.0924756614960671\e{-03}  &  2.3219312682036197\e{-02}  &  5.0546454457410162\e{-01} &  2.0474618401365854\e{+00} \\
 4.1043824968249148\e{-01}  &  8.3534647700054046\e{-02}  &  -8.8422252029506054\e{-01}  &  6.2852588972458059\e{-02}  &  5.4668509293072887\e{-01} &  -3.2336329115436924\e{-01} \\
 8.4898255952298962\e{-01}  &  9.8972579458252483\e{-02}  &  -8.9129367099545231\e{-01}  &  5.4473719351268962\e{-02}  &  7.1414182420995431\e{-01} &  3.2899060754742177\e{-01} \\
 3.3543896258348421\e{-01}  &  4.3010116145097040\e{-02}  &  1.5297157134040762\e{+00}  &  2.4345446089014514\e{-02}  &  -1.0558095282893749\e{+00} &  0.0000000000000000\e{+00} \\
\end{tabular}}
\caption{3S* low storage coefficients of ERK(20,5) method, optimized for 5$^\mathrm{th}$-order SD scheme.}
\label{tbl:rk-coeffs-3Sstar-erk5-20}
\end{table} 


\FloatBarrier

\newpage
\bibliography{optimized-erk-sd}

\end{document}